\documentclass[12pt]{amsart}
\usepackage{amsmath}
\usepackage{amssymb}
\usepackage{amsmath,amssymb,amsthm,amscd}
\usepackage{relsize}

\usepackage{txfonts}
\usepackage{amsbsy}
\usepackage{graphicx,color}

\numberwithin{equation}{section}

\newtheorem{prop}{Proposition}[section]
\newtheorem{theo}{Theorem}[section]
\newtheorem{lemm}{Lemma}[section]
\newtheorem{coro}{Corollary}[section]

\def\begeq{\begin{equation}}
\def\endeq{\end{equation}}

\begin{document}

\title{Bernstein Problem of Affine Maximal Type Hypersurfaces on Dimension $N\geq3$}
\author{Shi-Zhong Du}
\thanks{The author is partially supported by STU Scientific Research Foundation for Talents (SRFT-NTF16006), and Natural Science Foundation of Guangdong Province (2019A1515010605)}
  \address{The Department of Mathematics,
            Shantou University, Shantou, 515063, P. R. China.} \email{szdu@stu.edu.cn}

\renewcommand{\subjclassname}{%
  \textup{2010} Mathematics Subject Classification}
\subjclass[2010]{53A15 $\cdot$ 53A10 $\cdot$ 35J60}
\date{Mar. 2020}
\keywords{Affine maximal hypersurfaces, Bernstein theorem.}

\begin{abstract}
   Bernstein problem for affine maximal type equation
     \begin{equation}\label{e0.1}
      u^{ij}D_{ij}w=0, \ \ w\equiv[\det D^2u]^{-\theta},\ \ \forall x\in\Omega\subset{\mathbb{R}}^N
     \end{equation}
  has been a core problem in affine geometry. A conjecture proposed firstly by Chern (Proc. Japan-United States Sem., Tokyo, 1977, 17-30) for entire graph and then extended by Trudinger-Wang (Invent. Math., {\bf140}, 2000, 399-422) to its fully generality asserts that any Euclidean complete, affine maximal type, locally uniformly convex $C^4$-hypersurface in ${\mathbb{R}}^{N+1}$ must be an elliptic paraboloid. At the same time, this conjecture was solve completely by Trudinger-Wang for dimension $N=2$ and $\theta=3/4$, and later extended by Jia-Li (Results Math., {\bf56} 2009, 109-139) to $N=2, \theta\in(3/4,1]$ (see also Zhou (Calc. Var. PDEs., {\bf43} 2012, 25-44) for a different proof). On the past twenty years, much efforts were done toward higher dimensional issues but not really successful yet, even for the case of dimension $N=3$. In this paper, we will construct non-quadratic affine maximal type hypersurfaces which are Euclidean compete for
    $$
     N\geq3,\ \  \theta\in(1/2,(N-1)/N).
    $$
\end{abstract}

\maketitle\markboth{Affine Maximal Hypersurfaces}{Bernstein Theorem}

\tableofcontents

\section{Introduction}

In this paper, we study the local uniformly convex solution to affine maximal type equation
   \begin{equation}\label{e1.1}
       D_{ij}(U^{ij}w)=0, \ \ \forall x\in\Omega\subset{\mathbb{R}}^N,
   \end{equation}
where $U^{ij}$ is the co-factor matrix of $u_{ij}$ and $w\equiv[\det D^2u]^{-\theta}, \theta>0$. Equation \eqref{e1.1} is the Euler-Lagrange equation of the affine area functional
   \begin{eqnarray*}
     {\mathcal{A}}(u,\Omega)&\equiv&\int_\Omega[\det D^2u]^{1-\theta}\\
      &=&\int_{{\mathcal{M}}_\Omega}K_0^{1-\theta}(1+|Du|^2)^{\vartheta}dV_{g_0},\ \ \vartheta=\frac{N+1}{2}-\frac{N+2}{2}\theta
   \end{eqnarray*}
for $\theta\not=1$ and
    $$
     {\mathcal{A}}(u,\Omega)\equiv\int_\Omega\log\det D^2u
    $$
for $\theta=1$, where $g_0$ and $K_0$ are the induced metric and the Gauss curvature of the graph ${\mathcal{M}}_\Omega\equiv\Big\{(x,z)\in{\mathbb{R}}^{N+1}|\ z=u(x), x\in\Omega\Big\}$ respectively. Noting that
   $$
     D_jU^{ij}=0, \ \ \forall i=1,2,\cdots, N,
   $$
equation \eqref{e1.1} can also be written by
   $$
      U^{ij}D_{ij}w=0,
   $$
or equivalent to
   \begin{equation}\label{e1.2}
      u^{ij}D_{ij}w=0
   \end{equation}
for $[u^{ij}]$ denoting the inverse of usual metric $[u_{ij}]$ of graph ${\mathcal{M}}_\Omega$.

 The classical affine maximal case $\theta\equiv\frac{N+1}{N+2}$ has been studied extensively on the past. If one introduces the affine metric
   $$
     A_{ij}=\frac{u_{ij}}{[\det D^2u]^{1/(N+2)}}
   $$
on  ${\mathcal{M}}_\Omega$ comparing to the Calabi's metric $g_{ij}=u_{ij}$ and sets
  $$
    H\equiv[\det D^2u]^{-1/(N+2)},
  $$
it's not hard to see that \eqref{e1.2} turns to be
   \begin{equation}\label{e1.3}
    \triangle_{{\mathcal{M}}}H=0
   \end{equation}
for Laplace-Beltrami operator
   $$
     \triangle_{{\mathcal{M}}}\equiv\frac{1}{\sqrt{A}}D_i(\sqrt{A}A^{ij}D_j)=HD_i(H^{-2}u^{ij}D_j)
   $$
with respect to this affine metric, where $A$ is the determinant of $[A_{ij}]$ and $[A^{ij}]$ stands for the inverse of $[A_{ij}]$. So, the hypersurface ${\mathcal{M}}$ is affine maximal if and only if $H$ is harmonic on ${\mathcal{M}}_\Omega$.

A conjecture proposed by Chern \cite{Ch} for dimension $N=2$ and $\theta=3/4$ asserts that every locally convex entire graph of \eqref{e1.2} must be a paraboloid. Much efforts were done toward this conjecture, (see examples \cite{B,Ca1,Ca2,JL,TW} for partial results) until a landmark paper by Trudinger-Wang \cite{TW}. At there, they strengthened the Bernstein problem to its full generality  and then give a proof for 2 dimensional affine maximal case. Before our discussion, let's first reformulate this version of full Bernstein problem by Trudinger-Wang to all dimension $N$ and positive $\theta$ as following.\\

\noindent\textbf{Full Bernstein Problem in Sense of Trudinger-Wang:} Given dimension $N\geq1$ and $\theta>0$, whether any locally uniform convex, Euclidean complete affine maximal type hypersurfaces must be an elliptic paraboloid?\\

As shown in Theorem \ref{t2.1} in Section 2, for any given $\theta>0$, there exists a critical dimension $N^*(\theta)$ such that the Bernstein theorem holds for $1\leq N\leq N^*(\theta)$ and fails to hold for $N>N^*(\theta)$. Now, let's first recall a known Bernstein theorem on dimension $N\leq2$, which is mainly owe to Trudinger-Wang \cite{TW} for $N=2, \theta=3/4$. (see also \cite{JL2,Z} for $N=2, \theta\in(3/4,1]$)\\

  \noindent\textbf{Theorem A.} Full Bernstein theorem in sense of Trudinger-Wang holds under either one of the following cases:

  (1) $N=1$ and $\theta>0$, or

  (2) $N=2$ and $\theta\in[3/4,1]$.\\

Unlike the two dimensional case, attempts towards full Bernstein theorem on dimension $N\geq3$ were not yet successful. See examples \cite{A,D1,D2,D3,DF,JL,Z} for partial results in positive directions. In this paper, we will give an opposite answer to the Bernstein theorem for $N\geq3$ as following.\\

 \noindent\textbf{Theorem B.} For any $N\geq3$ and
   $$
    \theta\in(1/2,(N-1)/N),
   $$
  there exists a convex set $\Omega\subset{\mathbb{R}}^N$ and a non-quadratic $C^5$-solution $u$ of \eqref{e1.1} on $\Omega$ satisfying
    $$
     \lim_{x\to\partial\Omega}u(x)=+\infty.
    $$

  It's remarkable that a $W^{2,1}_{loc}({\mathbb{R}}^{10})$ solution has been known in \cite{TW} which violates the validity of Bernstein theorem, although it is in weak sense. On another hand, comparison of our result with some other known results under completeness of metric displays delicate differences between the metrics we choose. For examples, Jia-Li \cite{JL} and McCoy \cite{M} proved a Bernstein property for dimension $N=2,3$ and $\theta=(N+1)/(N+2)$ under completeness of Calabi's metric. Later, Zhou \cite{Z} extended the result to Abreu's equation $\theta=1$ for $2\leq N\leq4$. Recently, we \cite{DF} generalized the result of \cite{Z} to a wider range of $\theta$ for all dimension $N\geq2$, which covers the affine maximal equation for $N=2,3$ and the Abreu's equation for $2\leq N\leq5$.

    When restricting to the family of rotational symmetric solutions, we still have the following Bernstein theorem.\\

 \noindent\textbf{Theorem C.} For $N\geq3$ and $\theta>0$, the full Bernstein theorem holds for radial symmetric solutions.\\

This paper is organized as follows: Theorem A for $N=1$ will be proven in Proposition \ref{p2.1} of Section 2, and Theorem C will be proven in Theorem \ref{t3.1} of Section 3. Finally, we will prove Theorem B in Section 4-8.

\vspace{40pt}

\section{Critical dimension}

At first, let's verify the validity of full Bernstein theorem for dimension $N=1$ and all $\theta>0$.

\begin{prop}\label{p2.1}
  For $N=1$ and $\theta>0$, the full Bernstein theorem in Trudinger-Wang sense holds.
\end{prop}

\noindent\textbf{Proof.} We will first consider the case $\Omega={\mathbb{R}}$. For $N=1$, \eqref{e1.1} changes to O.D.E.
   \begin{equation}\label{e2.1}
      u^{(4)}=(\theta+1)\frac{(u''')^2}{u''}, \ \ \forall x\in{\mathbb{R}}.
   \end{equation}
 Near the location where $u'''$ is not identical to zero, we divided \eqref{e2.1} by $u'''$ and then integrate it in $x$. It gives that
   $$
    \log u'''-(\theta+1)\log u''=C_1
   $$
for constant $C_1$. Equivalently, we can obtain that
    $$
      u'''=C_2(u'')^{\theta+1}
    $$
and solve it by
   \begin{equation}\label{e2.2}
    (u'')^{-\theta}=-\theta C_2x+C_3.
   \end{equation}
Since $u$ is convex, the unique possibility is $C_2=0$, and thus the conclusion follows. In case $\Omega=[a,b], a<b\in{\mathbb{R}}$, we still have \eqref{e2.2} on $[a,b]$. Taking into account the large condition
   \begin{equation}\label{e2.3}
     \lim_{x\to\partial\Omega}u(x)=+\infty,
   \end{equation}
By translation invariant of \eqref{e1.1}, without loss of generality, we may assume that $a=0, b=1$. Then convexity of $u$ means that
  \begin{equation}\label{e2.4}
  \begin{cases}
    C_3\geq0, \ \ -\theta C_2+C_3>0\\
    \mbox{or } C_3>0, \ \ -\theta C_2+C_3\geq0.
  \end{cases}
  \end{equation}
However, condition \eqref{e2.4} will violate the validity of large condition \eqref{e2.3}, and thus exclude the possibility of $\Omega=[a,b]$ in this section.

 Finally, we consider the case $\Omega=[0,+\infty)$ after normalization. In this case, the convexity of $u$ implies that
  \begin{equation}\label{e2.5}
  \begin{cases}
    C_3\geq0, \ \ C_2<0\\
    \mbox{or } C_3>0, \ \ C_2\leq0.
  \end{cases}
  \end{equation}
 The large condition was violated also. So, the proposition holds true. $\Box$\\

The next theorem clarifies the existence of critical dimension for validity of Bernstein theorem.

\begin{theo}\label{t2.1}
  Suppose that $\varphi$ is a solution to \eqref{e1.1} on ${\mathbb{R}}^n$, then
     $$
      u(x,y)\equiv\varphi(x)+\frac{1}{2}|y|^2, \ \ y=(y^1,\cdots, y^m)
     $$
  is also a solution on ${\mathbb{R}}^{n+m}$, where $m\in{\mathbb{N}}$. As a corollary, for any $\theta>0$, there exists a critical dimension $N_*(\theta)\in{\mathbb{N}}\cup\{\infty\}$, such that Bernstein theorem holds for $1\leq N\leq N^*(\theta)$ and fails to hold for $N>N^*(\theta)$. Furthermore, if $\theta\in(1/2,1)$, we have
    \begin{equation}\label{e2.16}
      \Bigg[\frac{1}{1-\theta}\Bigg]\geq N_*(\theta)\geq\begin{cases}
         2, & \mbox{if } \frac{3}{4}\leq\theta<1,\\
         1, & \mbox{if } \frac{1}{2}<\theta<\frac{3}{4},
      \end{cases}
    \end{equation}
  where $[z]$ stands for the largest integer no greater than $z$.
\end{theo}

\noindent\textbf{Proof.} Noting that
   $$
    D^2u(x,y)=\left(
      \begin{array}{cc}
        D^2\varphi(x) & 0\\
        0 &  E_m
      \end{array}
      \right)
   $$
for unit matrix $E_m\equiv(\delta_{ij})_{m\times m}$ and

  $$
    w_u(x,y)=[\det D^2u]^{-\theta}=[\det D^2\varphi(x)]^{-\theta}\equiv w_\varphi(x),
  $$
we have
   $$
    u^{\alpha\beta}D_{\alpha\beta}w_u(x)=\varphi^{ij}D_{ij}w_\varphi(x)=0,\ \ i,j=1,\cdots,n, \ \ \alpha,\beta=1,\cdots,n+m.
   $$
Therefore, if there exists a non-quadratic solution to \eqref{e1.1} on ${\mathbb{R}}^n$, then there exists also non-quadratic solution to \eqref{e1.1} on ${\mathbb{R}}^{n+m}$ for any $m\in{\mathbb{N}}$. The proof of existence of critical dimension was done.

 Finally, \eqref{e2.16} is a direct consequence of Theorem A and B in introduction. $\Box$\\

\vspace{40pt}

\section{Group $SO(N)$-invariant solutions}

When considering solution of \eqref{e1.1} which is invariant under group action $SO(N)$, one has $u=u(|x|)$ and
   \begin{eqnarray*}
    D_iu&=&\frac{x_i}{|x|}u',\\
    D_{ij}u&=&\frac{x_ix_j}{|x|^2}u''+\Bigg(\frac{\delta_{ij}}{|x|}-\frac{x_ix_j}{|x|^3}\Bigg)u'\\
    &=&\frac{u'}{r}\Bigg\{\delta_{ij}+\Bigg(\frac{u''}{ru'}-\frac{1}{r^2}\Bigg)x_ix_j\Bigg\}.
   \end{eqnarray*}
Direct computation shows that
   \begin{eqnarray*}
    \det[D^2u]&=&\Bigg(\frac{u'}{r}\Bigg)^N\Bigg\{1+\Bigg(\frac{u''}{ru'}-\frac{1}{r^2}\Bigg)r^2\Bigg\}=u''\Bigg(\frac{u'}{r}\Bigg)^{N-1},\\
    u^{ij}&=&\Bigg(\frac{u'}{r}\Bigg)^{-1}\Bigg\{\delta_{ij}-\frac{ru''-u'}{r^3u''}x_ix_j\Bigg\},\\
    D_i&=&\frac{x_i}{r}\partial_r, \ \ D_{ij}=\frac{x_ix_j}{r^2}\partial^2_r+\Bigg(\frac{\delta_{ij}}{r}-\frac{x_ix_j}{r^3}\Bigg)\partial_r,
   \end{eqnarray*}
and thus
   \begin{eqnarray*}
     u^{ij}D_{ij}&=&\Bigg(\frac{u'}{r}\Bigg)^{-1}\Bigg(\partial^2_r+\frac{N-1}{r}\partial_r\Bigg)-\Bigg(\frac{u'}{r}\Bigg)^{-1}\frac{ru''-u'}{r^3u''}r^2\partial^2_r\\
      &=&\Bigg(\frac{u'}{r}\Bigg)^{-1}\Bigg\{\frac{u'}{ru''}\partial_r^2+\frac{N-1}{r}\partial_r\Bigg\}.
   \end{eqnarray*}
As a result, \eqref{e1.1} can be reformulated by
  \begin{eqnarray}\label{e3.1}
    && \ \ \ \ \ \ \ \ \ \  \ \ \ \ \ -u^{(4)}+(\theta+1)(u'')^{-1}(u''')^2+2(N-1)u'''\Bigg\{(\theta-1)\frac{u''}{u'}-\theta\frac{1}{r}\Bigg\}\\ \nonumber
     && \ \ \ \ \ \  +(N-1)u''\Bigg(\frac{u''}{u'}-\frac{1}{r}\Bigg)\Bigg\{\Big[(N-1)\theta-(N-2)\Big]\frac{u''}{u'}-\Big[(N-1)\theta-1\Big]\frac{1}{r}\Bigg\}=0
  \end{eqnarray}
Setting $v=u'$, \eqref{e3.1} changes to
   \begin{eqnarray}\label{e3.2}
     &\displaystyle -\frac{v'''}{v}+(\theta+1)\Bigg(\frac{v'}{v}\Bigg)^{-1}\Bigg(\frac{v''}{v}\Bigg)^2+2(N-1)\frac{v''}{v}\Bigg\{(\theta-1)\frac{v'}{v}-\theta\frac{1}{r}\Bigg\}&\\ \nonumber
     & \displaystyle +(N-1)\Bigg(\frac{v'}{v}\Bigg)\Bigg(\frac{v'}{v}-\frac{1}{r}\Bigg)\Bigg\{\Big[(N-1)\theta-(N-2)\Big]\frac{v'}{v}-\Big[(N-1)\theta-1]\frac{1}{r}\Bigg\}=0.&
   \end{eqnarray}
If one denotes $w=\frac{v'}{v}$, then
     $$
      \frac{v''}{v}=w'+w^2, \ \ \frac{v'''}{v}=w''+3ww'+w^3.
     $$
  Substituting into \eqref{e3.2}, it yields that
     \begin{eqnarray*}
       && -w''-3ww'-w^3+(\theta+1)w^{-1}(w'+w^2)^2+2(N-1)(w'+w^2)\Bigg\{(\theta-1)w-\theta\frac{1}{r}\Bigg\}\\
       &&+(N-1)w\Big(w-\frac{1}{r}\Big)\Bigg\{\Big[(N-1)\theta-(N-2)\Big]w-\Big[(N-1)\theta-1\Big]\frac{1}{r}\Bigg\}=0.
     \end{eqnarray*}
 Transforming $w(r)=\frac{1}{r}\eta(\log r)$, one has
     \begin{eqnarray*}
      w&=&\frac{1}{r}\eta\Leftrightarrow rw=\eta,\\
      w_r&=&-\frac{1}{r^2}\eta+\frac{1}{r^2}\eta'\Leftrightarrow r^2w_r=\eta'-\eta,\\
      w_{rr}&=&\frac{2}{r^3}\eta-\frac{3}{r^3}\eta'+\frac{1}{r^3}\eta''\Leftrightarrow r^3w_{rr}=\eta''-3\eta'+2\eta.
     \end{eqnarray*}
 Replacing $w$ by $\eta$, one concludes that
    \begin{eqnarray*}
      && -(\eta''-3\eta'+2\eta)-3\eta(\eta'-\eta)-\eta^3+(\theta+1)\eta^{-1}(\eta'-\eta+\eta^2)^2\\
      &&+2(N-1)(\eta'-\eta+\eta^2)\Bigg\{(\theta-1)\eta-\theta\Bigg\}\\
      &&+(N-1)\eta(\eta-1)\Bigg\{\Big[(N-1)\theta-(N-2)\Big]\eta-\Big[(N-1)\theta-1\Big]\Bigg\}=0
    \end{eqnarray*}
Letting $\zeta(\eta)=\eta'$, then
     $$
      \eta''=\frac{d\zeta}{d\eta}\eta'=\zeta\zeta'.
     $$
This will yield an equivalent form
    \begin{eqnarray}\nonumber\label{e3.3}
       && -\zeta\zeta'+(\theta+1)\eta^{-1}\zeta^2+\zeta\Bigg\{\Big[2N\theta-(2N-1)\Big]\eta-\Big[2N\theta-1\Big]\Bigg\}\\
       && \ \ \ \ \ \ \ +N\eta(\eta-1)\Bigg\{\Big[N\theta-(N-1)\Big]\eta-\Big[N\theta-1\Big]\Bigg\}=0
     \end{eqnarray}
of \eqref{e3.1}. We have the following necessary conditions subjected to strict convexity and smoothness of original solution $u$ of \eqref{e1.1}.\\

\begin{lemm}\label{l3.1}
  (1) The strict convexity of $u=u(|x|)$ is equivalent to $u_{rr}(r)>0$ holds for all $r\geq0$.

  (2) If $u=u(|x|)$ is a $C^2-$function, then
      $$
       u_r(0)=0\Leftrightarrow v(0)=0\Rightarrow \lim_{t\to-\infty}\eta(t)=1.
      $$
  Moreover, if $\zeta(\eta)=\eta'(t)>0$, we have
     $$
      \lim_{t\to-\infty}\eta'(t)=0\Rightarrow\lim_{\eta\to1^+}\zeta(\eta)=0.
     $$
  And if $\zeta(\eta)=\eta'(t)<0$, we have
     $$
      \lim_{t\to-\infty}\eta'(t)=0\Rightarrow\lim_{\eta\to1^-}\zeta(\eta)=0.
     $$
\end{lemm}

\noindent\textbf{Proof.} (1) Suppose that $u_{rr}(r)>0$ holds for all $r\geq0$, since $u_r(0)=0$ by smoothness of $u$ at origin, one has $u_r(r)>0$ for all $r>0$. Taking any $\xi\not=0\in{\mathbb{R}}^N$, there holds
   \begin{eqnarray}\nonumber\label{e3.4}
     u_{ij}\xi^i\xi^j&=&\frac{u'}{r}\Bigg\{\delta_{ij}+\Bigg(\frac{u''}{ru'}-\frac{1}{r^2}\Bigg)x_ix_j\Bigg\}\xi^i\xi^j\\
      &=&\frac{u'}{r}\Bigg\{|\xi|^2+\Bigg(\frac{u''}{ru'}-\frac{1}{r^2}\Bigg)(x\cdot\xi)^2\Bigg\},
   \end{eqnarray}
where $\frac{u'}{r}$ is understood as $u''$ at $r=0$ and
   $$
    \lim_{r\to0^+}r^2\Bigg(\frac{u''}{ru'}-\frac{1}{r^2}\Bigg)=0.
   $$
Now, we will focus on the case outside origin. The argument at origin is similar. If $\frac{u''}{ru'}-\frac{1}{r^2}\geq0$, it's clear that $u_{ij}\xi^i\xi^j>0$. If $\frac{u''}{ru'}-\frac{1}{r^2}<0$, then
   \begin{eqnarray*}
    u_{ij}\xi^i\xi^j&\geq&\frac{u'}{r}\Bigg\{|\xi|^2+\Bigg(\frac{u''}{ru'}-\frac{1}{r^2}\Bigg)r^2|\xi|^2\Bigg\}\\
      &=&\frac{u'}{r}\frac{ru''}{u'}|\xi|^2>0.
   \end{eqnarray*}
Conversely, if $u_{rr}(|x_0|)\leq0$ for some $0\not=x_0\in{\mathbb{R}}^N$, let's utilize \eqref{e3.4} for $\xi=x$. One deduces that $u_{ij}\xi^i\xi^j\Big|_{x_0}=0$, which contradicts with the assumption of strict convexity.

 To show (2), we need only use the Taylor expansion for smooth function together with the strict convexity of $u$. $\Box$\\

 \noindent\textbf{Remark 3.1} Above calculation shows that an equivalent form of Bernstein property for radial symmetric solution is given by the assertion that there is no non-trivial solution $\zeta$ of \eqref{e3.3} excepting the singular one
     $$
      \eta\equiv 1, \ \ \zeta\equiv0.
     $$

 \begin{theo}\label{t3.1}
    For dimension $N\geq3$ and all $\theta>0$, the unique solution of \eqref{e3.3} is given by the degenerate one
       \begin{equation}\label{e3.5}
        \zeta\equiv0, \ \ \eta\equiv1,
       \end{equation}
    which corresponds to solution $Cr^2, C>0$ of \eqref{e3.1}.
 \end{theo}

\noindent\textbf{Proof.} The solutions of \eqref{e3.3} can be divided into three types. The first one is the degenerate solution \eqref{e3.5}. The second one is the positive solution on $(\eta_k,\eta_k')$ for
  $$
   \eta_k\geq1, \ \ \lim_{k\to+\infty}\eta_k=1,
  $$
which satisfies that
   $$
    \zeta(\eta_k)=0, \ \ \zeta(\eta)>0, \ \ \forall \eta\in(\eta_k,\eta'_k).
   $$
 The last one is the negative solution on $(\eta'_k,\eta_k)$ for
  $$
   \eta_k\leq1, \ \ \lim_{k\to+\infty}\eta_k=1,
  $$
which satisfies that
   $$
    \zeta(\eta_k)=0, \ \ \zeta(\eta)<0, \ \ \forall \eta\in(\eta'_k,\eta_k).
   $$
To exclude the possibility of type 2, we set $\varphi(\eta)\equiv\eta^{-2(\theta+1)}\zeta^2(\eta), \zeta>0$. Then \eqref{e3.3} changes into
  \begin{eqnarray}\label{e3.6}
   &-\varphi'+2\sqrt{\varphi}\eta^{-(\theta+1)}\Bigg\{\Big[2N\theta-(2N-1)\Big]\eta-\Big[2N\theta-1\Big]\Bigg\}&\\ \nonumber
   & +2N\eta^{-(2\theta+1)}(\eta-1)\Bigg\{\Big[N\theta-(N-1)\Big]\eta-\Big[N\theta-1\Big]\Bigg\}=0, \ \ \forall \eta\in(\eta_k,\eta'_k).&
  \end{eqnarray}
Noting that for $N\geq3$ and $\theta>0$, it's inferred from \eqref{e3.6} that
   $$
    \varphi(\eta_k)=0, \ \ \varphi'(\eta)<0, \ \ \forall\eta\in(\eta_k,\eta'_k)
   $$
for $k$ large. This contradicts with the positivity of $\varphi$ on $(\eta_k,\eta'_k)$. Similarly, to exclude the possibility of type 3, we set $\varphi(\eta)\equiv\eta^{-2(\theta+1)}\zeta^2(\eta), \zeta<0$. Then \eqref{e3.3} changes into
  \begin{eqnarray}\label{e3.7}
   &-\varphi'-2\sqrt{\varphi}\eta^{-(\theta+1)}\Bigg\{\Big[2N\theta-(2N-1)\Big]\eta-\Big[2N\theta-1\Big]\Bigg\}&\\ \nonumber
   & +2N\eta^{-(2\theta+1)}(\eta-1)\Bigg\{\Big[N\theta-(N-1)\Big]\eta-\Big[N\theta-1\Big]\Bigg\}=0, \ \ \forall \eta\in(\eta_k,\eta'_k).&
  \end{eqnarray}
Noting that for $N\geq3$ and $\theta>0$, it's inferred from \eqref{e3.7} that
   $$
    \varphi(\eta_k)=0, \ \ \varphi'(\eta)>0, \ \ \forall\eta\in(\eta'_k,\eta_k)
   $$
for $k$ large. This also contradicts with the positivity of $\varphi$ on $(\eta'_k,\eta_k)$. So, the uniqueness of \eqref{e3.5} has been proven. $\Box$\\

\begin{coro}\label{c3.1}
  For $N=2k, k\geq2, k\in{\mathbb{Z}}$ and $\theta=\frac{N+1}{N+2}$,
    \begin{equation}\label{e3.8}
      u(x)=C|x|^{2k^2}, \ \ \forall x\in{\mathbb{R}}^N
    \end{equation}
  are all entire solutions to \eqref{e1.1} on ${\mathbb{R}}^N\setminus\{0\}$.\\
\end{coro}

\noindent\textbf{Proof.} The proof is easy by noting that $\eta\equiv\frac{N\theta-1}{N\theta-(N-1)}$, $\zeta=\eta'\equiv0$ is a degenerate solution to \eqref{e3.3}. Therefore, \eqref{e3.8} yields smooth solution to \eqref{e1.1} outside origin. $\Box$\\

\noindent\textbf{Remark 3.2} The solution in Corollary \ref{c3.1} does not produce a true counter example of Bernstein problem since the convexity of $u$ is degenerate at $x=0$.\\

\noindent\textbf{Remark 3.3} Combining Corollary \ref{c3.1} with Theorem \ref{t2.1}, one actually obtains examples for all dimension $N\geq4$.

\vspace{40pt}

\section{Reduction to opposite pairs of nonlinear eigenvalue problem}

In this section, we attempt to construct true counter examples of \eqref{e1.1}. Suppose $x\in{\mathbb{R}}^n$ and $y\in{\mathbb{R}}^m$, we look for solution in form of $u(x,y)=\varphi(x)+\psi(y)$ on ${\mathbb{R}}^N, N=n+m$. Setting $w_\varphi\equiv[\det D^2\varphi(x)]^{-\theta}$ and $ w_\psi\equiv[\det D^2\psi(y)]^{-\theta}$ respectively, we have
   $$
    D^2u=\left(
     \begin{array}{cc}
       D^2\varphi(x) & 0 \\
       0 & D^2\psi(y)
     \end{array}
     \right)
   $$
and so \eqref{e1.1} is changed to
   \begin{equation}\label{e4.1}
    \Big(\varphi^{ij}(x)D_{ij}w_\varphi(x)\Big)w_\psi(y)+\Big(\psi^{ij}(y)D_{ij}w_\psi(y)\Big)w_\varphi(x)=0.
   \end{equation}
Therefore, there must be a nonnegative constant $\lambda$ such that
  \begin{equation}\label{e4.2}
   \frac{\varphi^{ij}(x)D_{ij}w_\varphi(x)}{w_\varphi(x)}=-\frac{\psi^{ij}(y)D_{ij}w_\psi(y)}{w_\psi(y)}=-\lambda.
  \end{equation}
So, finding of non-quadratic solution of \eqref{e1.1} is reduced to looking for non-trivial solutions to nonlinear eigenvalue problem
  \begin{equation}\label{e4.3}
    u^{ij}D_{ij}w=\lambda w, \ \ \forall x\in{\mathbb{R}}^n
  \end{equation}
for opposite pairs $(\lambda, n)$ and $(-\lambda, m)$,
where $n,m$ are positive integers and $w\equiv[\det D^2u]^{-\theta}$.\\

\noindent\textbf{Remark 4.1} Noting that the constant $\lambda$ in \eqref{e4.3} can be adjusted when one scales $u\to\kappa u$, so false of Bernstein theorem is reduced to find non-trivial solutions of \eqref{e4.3} for some positive pair of $(\lambda_1,n), \lambda_1>0$ and negative pair $(\lambda_2,m), \lambda_2<0$, where $n, m$ are positive integers.\\

 Supposing $u=u(|x|)$, \eqref{e4.3} can be reformulated to
  \begin{eqnarray}\label{e4.4}
    &-u^{(4)}+(\theta+1)(u'')^{-1}(u''')^2+2(n-1)u'''\Big\{(\theta-1)\frac{u''}{u'}-\theta\frac{1}{r}\Big\}&\\ \nonumber
     &+(n-1)u''\Big(\frac{u''}{u'}-\frac{1}{r}\Big)\Bigg\{\Big[(n-1)\theta-(n-2)\Big]\frac{u''}{u'}-\Big[(n-1)\theta-1\Big]\frac{1}{r}\Bigg\}=\lambda'(u'')^2,&
  \end{eqnarray}
as in Section 3, where $\lambda'\equiv-\frac{\lambda}{\theta}$. Setting $v=u'$, we change \eqref{e4.4} into
   \begin{eqnarray}\label{e4.5}
     &\displaystyle -\frac{v'''}{v}+(\theta+1)\Bigg(\frac{v'}{v}\Bigg)^{-1}\Bigg(\frac{v''}{v}\Bigg)^2+2(n-1)\frac{v''}{v}\Bigg\{(\theta-1)\frac{v'}{v}-\theta\frac{1}{r}\Bigg\}&\\ \nonumber
     & \displaystyle +(n-1)\Bigg(\frac{v'}{v}\Bigg)\Bigg(\frac{v'}{v}-\frac{1}{r}\Bigg)\Bigg\{\Big[(n-1)\theta-(n-2)\Big]\frac{v'}{v}-\Big[(n-1)\theta-1]\frac{1}{r}\Bigg\}=\lambda'\frac{(v')^2}{v}.&
   \end{eqnarray}
Denoting $w=\frac{v'}{v}$, simple calculation shows that
     $$
      \frac{v''}{v}=w'+w^2, \ \ \frac{v'''}{v}=w''+3ww'+w^3.
     $$
Substituting into \eqref{e4.5}, it yields that
     \begin{eqnarray*}
       && -w''-3ww'-w^3+(\theta+1)w^{-1}(w'+w^2)^2+2(n-1)(w'+w^2)\Bigg\{(\theta-1)w-\theta\frac{1}{r}\Bigg\}\\
       &&+(n-1)w\Big(w-\frac{1}{r}\Big)\Bigg\{\Big[(n-1)\theta-(n-2)\Big]w-\Big[(n-1)\theta-1\Big]\frac{1}{r}\Bigg\}=\lambda''w^2 e^{\int_1^r w(r)dr}.
     \end{eqnarray*}
After transforming $w(r)=\frac{1}{r}\eta(\log r)$, one gets that
     \begin{eqnarray*}
      w&=&\frac{1}{r}\eta\Leftrightarrow rw=\eta,\\
      w_r&=&-\frac{1}{r^2}\eta+\frac{1}{r^2}\eta'\Leftrightarrow r^2w_r=\eta'-\eta,\\
      w_{rr}&=&\frac{2}{r^3}\eta-\frac{3}{r^3}\eta'+\frac{1}{r^3}\eta''\Leftrightarrow r^3w_{rr}=\eta''-3\eta'+2\eta.
     \end{eqnarray*}
Replacing $w$ by $\eta$ in above formula again, one concludes that
    \begin{eqnarray}\nonumber\label{a4.6}
      & -(\eta''-3\eta'+2\eta)-3\eta(\eta'-\eta)-\eta^3+(\theta+1)\eta^{-1}(\eta'-\eta+\eta^2)^2&\\
      &+(n-1)\eta(\eta-1)\Bigg\{\Big[(n-1)\theta-(n-2)\Big]\eta-\Big[(n-1)\theta-1\Big]\Bigg\}&\\ \nonumber
      &+2(n-1)(\eta'-\eta+\eta^2)\Bigg\{(\theta-1)\eta-\theta\Bigg\}=\lambda''\eta^2e^{t+\int^t_0 \eta(t)dt}.&
    \end{eqnarray}
Letting $\zeta(\eta)=\eta'$, then
     $$
      \eta''=\frac{d\zeta}{d\eta}\eta'=\zeta\zeta'.
     $$
This will yield an equivalent form
    \begin{eqnarray}\nonumber\label{e4.7}
       && -\zeta\zeta'+(\theta+1)\eta^{-1}\zeta^2+\zeta\Bigg\{\Big[2n\theta-(2n-1)\Big]\eta-\Big[2n\theta-1\Big]\Bigg\}\\
       && \ \ \ \ \ \ \ +n\eta(\eta-1)\Bigg\{\Big[n\theta-(n-1)\Big]\eta-\Big[n\theta-1\Big]\Bigg\}=\lambda'''\eta^2\exp{\int^\eta_{\eta_0}\frac{s+1}{\zeta(s)}ds}
     \end{eqnarray}
of \eqref{e4.4} for $\eta, \eta_0>1$ in case $\zeta>0$ and $\eta, \eta_0<1$ in case $\zeta<0$. Next proposition constructs non-quadratic solution of \eqref{e4.3} for positive pair $(\lambda,1)$.

\begin{prop}\label{p4.1}
  For $n=1$ and $\theta>1/2$, there exists a non-quadratic solution $u\in C^\infty({\mathbb{R}})$ to \eqref{e4.3} for positive pair $(\lambda,1), \lambda>0$.
\end{prop}

\noindent\textbf{Proof.} When $n=1$, \eqref{e4.4} changes to
    $$
     -u^{(4)}+(\theta+1)\frac{(u''')^2}{u''}=\lambda(u'')^2.
    $$
 Setting $v=u''>0$, we get
   \begin{equation}\label{e4.8}
    -v''+(\theta+1)\frac{(v')^2}{v}=\lambda v^2.
   \end{equation}
 Regarding $v$ as variable and $w=v'$ as function, one has
   $$
    v''=\frac{\partial w}{\partial r}=\frac{\partial w}{\partial v}\frac{\partial v}{\partial r}=ww'.
   $$
 Therefore,
   \begin{eqnarray*}
     &&-ww'+(\theta+1)\frac{w^2}{v}=\lambda v^2\\
     &\Leftrightarrow&(w^2)'-2(\theta+1)\frac{w^2}{v}=-2\lambda v^2\\
     &\Leftrightarrow& w^2=\frac{2\lambda}{2\theta-1}v^3+C_0v^{2(\theta+1)}\equiv av^3+C_0v^{2(\theta+1)}, \ \ a\equiv\frac{2\lambda}{2\theta-1}.
   \end{eqnarray*}
Taking $C_0=-av_0^{1-2\theta}<0$ for positive initial datum $v_0$ of $v$ and noting that $w=v'<0$, we get
  \begin{equation}\label{e4.9}
   \begin{cases}
     v'=-\sqrt{a\big(v^3-v_0^{1-2\theta}v^{2(\theta+1)}\big)}, & \forall r\geq0\\
     v(0)=v_0.
   \end{cases}
  \end{equation}
Since $\theta>1/2$, this yields an unique positive solution $v\in(0,v_0)$ of \eqref{e4.9} by
  \begin{equation}\label{e4.10}
    \int_{v(r)}^{v_0}\frac{dv}{\sqrt{v^3-v_0^{1-2\theta}v^{2(\theta+1)}}}=\sqrt{a}r, \ \ \forall r\geq0,
  \end{equation}
thanks to
   $$
    \int_{0}^{v_0}\frac{dv}{\sqrt{v^3-v_0^{1-2\theta}v^{2(\theta+1)}}}=+\infty.
   $$
Setting
  $$
    u(r)=\int^r_0\int^t_0v(s)dsdt, \ \ \forall r\geq0,
  $$
it's clear that $u'(0)=0$ and $u'''(0)=v'(0)=0$. Taking more derivatives, we have
  \begin{eqnarray*}
    v''&=&-\frac{\sqrt{a}}{2}\frac{3v^2-2(\theta+1)v_0^{1-2\theta}v^{2\theta+1}}{\sqrt{v^3-v_0^{1-2\theta}v^{2(\theta+1)}}}v'\\
     &=&\frac{3a}{2}v^2-(\theta+1)av_0^{1-2\theta}v^{2\theta+1}
  \end{eqnarray*}
and
   \begin{eqnarray*}
    v'''&=&a\Big[3v-(\theta+1)(2\theta+1)v_0^{1-2\theta}v^{2\theta}\Big]v'\\
    v^{(4)}&=&a\Bigg\{\Big[3v-(\theta+1)(2\theta+1)v_0^{1-2\theta}v^{2\theta}\Big]v''+\Big[3-2\theta(\theta+1)(2\theta+1)v_0^{1-2\theta}v^{2\theta-1}\Big](v')^2\Bigg\}\\
    v^{(5)}&=&a\Bigg\{\Big[3v-(\theta+1)(2\theta+1)v_0^{1-2\theta}v^{2\theta}\Big]v'''+3\Big[3-2\theta(\theta+1)(2\theta+1)v_0^{1-2\theta}v^{2\theta-1}\Big]v'v''\\
    &&-(2\theta-1)2\theta(\theta+1)(2\theta+1)v_0^{1-2\theta}v^{2\theta-2}(v')^3\Bigg\}.
   \end{eqnarray*}
By mathematical induction, it's not hard to verify that $v\in C^\infty([0,+\infty))$ and
   \begin{equation}\label{e4.11}
    u^{(2k-1)}(0)=0, \ \ \forall k\in{\mathbb{N}}.
   \end{equation}
Thus,  $u(|x|)$ is a smooth function on ${\mathbb{R}}$ satisfying \eqref{e4.3} for positive pair $(\lambda,1)$.\\

\noindent\textbf{Remark 4.2} Using a similar argument, it's not difficult to verify that when $\theta>1/2$, the unique solution $u$ of \eqref{e4.3} for negative pair $(-\lambda,1)$ must be given by
   $$
    \int^{v(r)}_{v_0}\frac{dv}{\sqrt{a\Big(v_0^{1-2\theta}v^{2(\theta+1)}-v^3\Big)}}=r, \ \ \forall r\geq0.
   $$
The solution exists until
  $$
   r=R\equiv\int^{+\infty}_{v_0}\frac{dv}{\sqrt{a\Big(v_0^{1-2\theta}v^{2(\theta+1)}-v^3\Big)}}<+\infty.
  $$
However, since $\theta>1/2$,
  $$
   u(r)=\int^r_0\int^s_0v(s)ds
  $$
is not a large function near $r=R$. So, we will turn to look for solutions of \eqref{e4.7} for negative pair $(-\lambda,n), n\geq2$ in next three sections.

\vspace{40pt}

\section{Compatible conditions of \eqref{e4.7} at origin}

 For radial symmetric convex solution $u(|x|)\in C^4({\mathbb{R}}^n)$ of \eqref{e4.3}, one needs at least
   \begin{eqnarray}\nonumber\label{e5.1}
    &u'(0)=u'''(0)=0, u''(0)>0\Leftrightarrow v(0)=v''(0)=0, v'(0)>0\\
    &\Rightarrow\eta|_{t=-\infty}=1, \ \ \eta'-\eta+\eta^2|_{t=-\infty}=0\Rightarrow \eta(-\infty)=1, \ \ \zeta(1)=0.&
   \end{eqnarray}
For sufficiency of \eqref{e5.1}, we have the following result.

\begin{prop}\label{p5.1}
Letting $\zeta\in C^1([1,\eta_0])$ be a solution to \eqref{e4.7} satisfying first compatible condition
    \begin{equation}\label{e5.2}
     \begin{cases}
      \zeta(1)=0, \ \ \zeta'(1)=2, \ \ \zeta(\eta)>0, \ \ \forall \eta\in(1,\eta_0],\\
      \zeta'(\eta) \mbox{ monotone non-decreasing in } [1,\eta_0]
     \end{cases}
    \end{equation}
  for some $\eta_0>1$, the recovery solution $u(|x|)\in C^3(B_1)$ satisfies that $Du(0)=D^3u(0)=0$.
\end{prop}

\noindent\textbf{Proof.} At first, $\eta(t)$ is given by
   \begin{equation}\label{e5.3}
     t=\int^\eta_{\eta_0}\frac{1}{\zeta(s)}ds, \ \ \forall \eta>1
   \end{equation}
 for $\eta(0)=\eta_0$. Secondly, since
    \begin{eqnarray*}
      &&(\log v)_r=\frac{\eta(\log r)}{r}\\
      &\Leftrightarrow&\log v-\log v_0=\int^r_{r_0}\frac{\eta(\log\tau)}{\tau}d\tau\\
      &\Leftrightarrow& v=v_0\exp\Bigg\{\int^{r}_{r_0}\frac{\eta(\log\tau)}{\tau}d\tau\Bigg\}=v_0\exp\Bigg\{\int^{\log r}_{\log r_0}\eta(t)dt\Bigg\}.
    \end{eqnarray*}
$u$ can be restored by
   \begin{equation}\label{e5.4}
     u(r)=v_0\int^r_0\exp\Bigg\{\int^{\log s}_{\log r_0}\eta(t)dt\Bigg\}ds+u_0, \ \ \forall r\in[0,1],
   \end{equation}
 where $\eta(t)$ is given by \eqref{e5.3}. Setting $\overline{\eta}(r)=\eta(\log r)$ for simplicity, then
   \begin{equation}\label{e5.5}
     v(r)=v_0\exp\Bigg\{\int^r_{r_0}\frac{\overline{\eta}(s)}{s}dr\Bigg\}, \ \ \forall r\in[0,1],
   \end{equation}
 where $\overline{\eta}$ is a solution to
   \begin{equation}\label{e5.6}
    \begin{cases}
     \displaystyle \frac{d\overline{\eta}}{dr}=\frac{\zeta(\overline{\eta})}{r}, & r\in[0,1]\\
      \overline{\eta}(0)=1.
    \end{cases}
   \end{equation}

 \vspace{10pt}

\noindent\textbf{Claim:} Under assumption of Proposition \ref{p5.1}, the solution $v(r)$ given by \eqref{e5.5} and \eqref{e5.6} belongs to $C^2([0,1])$ and satisfies that
   \begin{equation}\label{e5.7}
     v(0)=0, \ \ v'(0)>0, \ \ v''(0)=0.
   \end{equation}

\noindent\textbf{Proof.} Before prove the claim, let's prepare two lemmas.

\begin{lemm}\label{l5.1}
  Letting $\zeta$ be given as in Proposition \ref{p5.1}, there exists an unique non-trivial solution $\overline{\eta}\in C^2((0,1])\cap C^1([0,1])$ satisfying
    \begin{equation}\label{e5.8}
      \overline{\eta}(r)-1\leq C_\alpha r^{1+\alpha}, \ \ \forall r\in[0,1]
    \end{equation}
  for any $\alpha<1$, where $C_\alpha$ are positive constants depending on $\alpha$.
\end{lemm}

\noindent\textbf{Proof.} Given any $r_0>0$ and $\overline{\eta}_0>1$, the local solution $\overline{\eta}\in C^2((0,1])\cap C([0,1])$ of \eqref{e5.6} is given by
  \begin{equation}\label{e5.9}
    \int^{\overline{\eta}}_{\overline{\eta}_0}\frac{d\overline{\eta}}{\zeta(\overline{\eta})}=\log\frac{r}{r_0}, \ \ \forall r\in[0,1].
  \end{equation}
By \eqref{e5.2} and Lagrange's theorem,
   \begin{equation}\label{e5.10}
     \lim_{\eta\to1^+}\frac{\zeta(\overline{\eta})}{\overline{\eta}-1}=\lim_{\eta\to1^+}\zeta'(\xi)=2.
   \end{equation}
Therefore, if $r$ is small, we derive from \eqref{e5.9} that
  \begin{eqnarray*}
    \log\frac{r}{r_0}\geq\int^{\overline{\eta}}_{\overline{\eta}_0}\frac{d\overline{\eta}}{2\alpha(\overline{\eta}-1)}=\frac{1}{2\alpha}\Big(\log(\overline{\eta}-1)-\log(\overline{\eta}_0-1)\Big), \ \ \forall r\in(0,1],
  \end{eqnarray*}
where $\alpha$ is a positive constant which can be taken closing to $1$ arbitrarily from below as long as $r$ is small. Thus, one gets
   \begin{equation}\label{e5.11}
     \overline{\eta}(r)-1\leq C_{\alpha, r_0}r^{2\alpha}, \ \ \forall r\in[0,1].
   \end{equation}
Next, we show that $\overline{\eta}\in C^1([0,1])$. In fact, by \eqref{e5.11}, $\overline{\eta}'(0)=0$. Another hand, it's inferred from \eqref{e5.6} that
    $$
     \frac{d\overline{\eta}}{dr}=\frac{\zeta(\overline{\eta}(r))-\zeta(\overline{\eta}(0))}{r}=\zeta'(\xi)\frac{\overline{\eta}(r)-\overline{\eta}(0)}{r}\to0=\frac{d\overline{\eta}}{dr}\Big|_{r=0}
    $$
as $r\to0^+$, where $\xi\in(\overline{\eta}(0),\overline{\eta}(r))$. So, $\overline{\eta}\in C^1([0,1])$. $\Box$\\

\begin{lemm}\label{l5.2}
  The non-trivial solution $\overline{\eta}$ given in Lemma \ref{l5.1} belongs to $C^2([0,1])$ and satisfies that $\overline{\eta}'(0)=0$ and
    \begin{equation}\label{e5.12}
      0\leq1-\overline{\eta}(r)\leq Cr^2, \ \ \forall r\in[0,1]
    \end{equation}
  for some positive constant $C$.
\end{lemm}

\noindent\textbf{Proof.} Differentiating \eqref{e5.6} at $r$, one gets
  \begin{eqnarray*}
    \frac{d^2\overline{\eta}}{dr^2}&=&\frac{\zeta'(\overline{\eta})}{r}\frac{d\overline{\eta}}{dr}-\frac{\zeta(\overline{\eta})}{r^2}\\
     &=&\frac{\zeta'(\overline{\eta})-1}{r}\frac{d\overline{\eta}}{dr}\\
     &\geq&\frac{1}{r}\frac{d\overline{\eta}}{dr}, \ \ \forall r\in(0,r_0)
  \end{eqnarray*}
by convexity of $\zeta$ in $(1,\overline{\eta}(r_0))$ and initial condition $\zeta'(1)=2$, where $r_0$ is chosen small. Thus,
   $$
    \frac{d}{dr}\Bigg(r^{-1}\frac{d\overline{\eta}}{dr}\Bigg)\geq0\Leftrightarrow r^{-1}\frac{d\overline{\eta}}{dr}(r)\downarrow\mbox{ as } r\downarrow, \ \ \forall r\in(0,r_0).
   $$
As a result, we conclude that
   \begin{eqnarray*}
     && r^{-1}\frac{d\overline{\eta}}{dr}(r)\leq\frac{d\overline{\eta}}{dr}(1)=\zeta(\overline{\eta}(1))\\
     &\Rightarrow&\frac{d\overline{\eta}}{dr}(r)\leq Cr\\
     &\Rightarrow&\frac{\overline{\eta}(r)-1}{r}=\frac{d\overline{\eta}}{dr}(\xi)\leq Cr, \ \ \xi\in(0,r)\\
     &\Rightarrow& 0<\overline{\eta}(r)-1<Cr^2, \ \ \forall 0<r<r_0.
   \end{eqnarray*}
So, \eqref{e5.12} holds true. To show $\overline{\eta}\in C^2([0,1])$, noting first that by monotonicity and boundedness of
  $$
   \frac{\frac{d\overline{\eta}}{dr}(r)-\frac{d\overline{\eta}}{dr}(0)}{r}=\frac{1}{r}\frac{d\overline{\eta}}{dr}
  $$
from above, the second derivative $\frac{d^2\overline{\eta}}{dr^2}(0)$ exists. Another hand, using the relation
   $$
    \frac{d^2\overline{\eta}}{dr^2}=\frac{\zeta'(\overline{\eta})-1}{r}\frac{d\overline{\eta}}{dr}
   $$
obtained above, monotonicity of $r^{-1}\frac{d\overline{\eta}}{dr}(r)$ and $\overline{\eta}$, it's inferred from the convexity of $\zeta$ that $\frac{d^2\overline{\eta}}{dr^2}(r)$ is monotone non-increasing and bounded from below by zero. So, we derive that
   $$
    \lim_{r\to0^+}\frac{d^2\overline{\eta}}{dr^2}(r)=\frac{d^2\overline{\eta}}{dr^2}(0)
   $$
and hence $\overline{\eta}\in C^2([0,1])$. $\Box$\\

Now, we can complete the proof of our claim. In fact, by Lemma \ref{l5.2} and
   $$
    \overline{\eta}(r)=1+\int^r_0(r-s)\overline{\eta}''(s)ds,
   $$
we have
   \begin{eqnarray*}
     v(r)&=&v_0\exp\Bigg\{\int^r_{r_0}\frac{1+\int^s_0(s-\tau)\overline{\eta}''(\tau)d\tau}{s}ds\Bigg\}\\
      &=&\frac{v_0r}{r_0}\exp\Bigg\{\int^r_{r_0}\frac{(\overline{\eta}''(0)+o(s))\frac{s^2}{2}}{s}ds\Bigg\}=\frac{v_0r}{r_0}\exp\Bigg\{\int^r_{r_0}\frac{s}{2}(\overline{\eta}''(0)+o(s))ds\Bigg\}\\
     v'(r)&=&\frac{v_0}{r_0}\exp\Bigg\{\int^r_{r_0}\frac{s}{2}(\overline{\eta}''(0)+o(s))ds\Bigg\}\Bigg\{1+\frac{r}{2}\Big(\overline{\eta}''(0)+o(r)\Big)\Bigg\}\\
     v''(r)&=&\frac{v_0}{r_0}\exp\Bigg\{\int^r_{r_0}\frac{s}{2}(\overline{\eta}''(0)+o(s))ds\Bigg\}\Bigg\{\Big(\overline{\eta}''(0)+o(r)\Big)+\frac{ro'(r)}{2}+\frac{r^2}{4}\Big(\overline{\eta}''(0)+o(r)\Big)^2\Bigg\}
   \end{eqnarray*}
are both continuous functions satisfying
   $$
    \lim_{r\to0^+}v(r)=v(0)=0, \ \ \lim_{r\to0^+}v'(r)=v'(0)>0, \ \ \lim_{r\to0^+}v''(r)=v''(0)=0,
   $$
where
  $$
   o(r)\equiv\frac{\int^r_0(r-\tau)\overline{\eta}''(\tau)d\tau}{\int^r_0(r-\tau)d\tau}-\overline{\eta}''(0)
  $$
is a smooth function on $(0,+\infty)$ satisfying
  $$
   \lim_{r\to0^+}o(r)=\lim_{r\to0^+}ro'(r)=0.
  $$
The conclusion was drawn. $\Box$\\

So, Proposition \ref{p5.1} is a direct consequence of Claim 1. $\Box$\\

Next, we will enhance the regularity of $u(|x|)$ in the following proposition.

\begin{prop}\label{p5.2}
  Letting $\zeta\in C^3([1,\eta_0])\cap C^\infty((1,\eta_0])$ be a solution to \eqref{e4.7} satisfying second compatible condition
    \begin{equation}\label{e5.13}
     \begin{cases}
      \zeta(1)=0, \ \ \zeta'(1)=2,\\
     \zeta(\eta)>0, \ \  \zeta''(\eta)0,\ \ \zeta^{(4)}(\eta)>0,\ \  \forall \eta\in(1,\eta_0],
     \end{cases}
    \end{equation}
 we have $v\in C^4([0,1])$ satisfying
   $$
    v(0)=v''(0)=v^{(4)}(0)=0
   $$
 and hence $u(|x|)\in C^5(B_1)$.
\end{prop}

\noindent\textbf{Proof.} As in proving of Lemma \ref{l5.2}, let's take third derivative on $\overline{\eta}$ by
  \begin{eqnarray*}
    \frac{d^3\overline{\eta}}{dr^3}&=&\frac{\zeta'(\overline{\eta})-1}{r}\frac{d^2\overline{\eta}}{dr^2}+\Bigg[\frac{\zeta''(\overline{\eta})\overline{\eta}'}{r}-\frac{\zeta'(\overline{\eta})-1}{r^2}\Bigg]\frac{d\overline{\eta}}{dr}\\
     &=&\Bigg[\frac{\zeta'(\overline{\eta})-1}{r}+\frac{\zeta(\overline{\eta})\zeta''(\overline{\eta})-(\zeta'(\overline{\eta})-1)}{r(\zeta'(\overline{\eta})-1)}\Bigg]\frac{d^2\overline{\eta}}{dr^2}\\
     &=&\frac{\zeta(\overline{\eta})\zeta''(\overline{\eta})+(\zeta'(\overline{\eta})-1)(\zeta'(\overline{\eta})-2)}{r(\zeta'(\overline{\eta})-1)}\frac{d^2\overline{\eta}}{dr^2}\to0
  \end{eqnarray*}
as $r\to0^+$, where $\xi\in(0,r), \chi\in(0,\xi)$ and the convergence of
   \begin{eqnarray}\nonumber\label{e5.14}
     &\frac{\zeta(\overline{\eta})}{r^2/2}=\frac{\zeta'(\overline{\eta}(\xi))\overline{\eta}'(\xi)}{\xi}=\zeta''(\overline{\eta}(\chi))\overline{\eta}'(\chi)\overline{\eta}'(\xi)+\zeta'(\overline{\eta}(\xi))\overline{\eta}''(\chi)&\\
     &\frac{\zeta'(\overline{\eta})-2}{r^2/2}=\frac{\zeta''(\overline{\eta}(\xi))\overline{\eta}'(\xi)}{\xi}=\zeta''(\overline{\eta}(\xi))\overline{\eta}''(\chi)&
   \end{eqnarray}
as $r\to0$ has been used. As a result, we get $\overline{\eta}\in C^3([0,1])$ with $\overline{\eta}'''(0)=0$. Similarly,
  \begin{eqnarray*}
   &\displaystyle\frac{d^4\overline{\eta}}{dr^4}=\frac{\zeta(\overline{\eta})\zeta''(\overline{\eta})+(\zeta'(\overline{\eta})-1)(\zeta'(\overline{\eta})-2)}{r(\zeta'(\overline{\eta})-1)}\frac{d^3\overline{\eta}}{dr^3}+\frac{\zeta(\overline{\eta})\zeta''(\overline{\eta})-(\zeta'(\overline{\eta})-2)}{r^2}\frac{d^2\overline{\eta}}{dr^2}&\\
    &\displaystyle+\Bigg\{\frac{\zeta(\overline{\eta})[\zeta'(\overline{\eta})\zeta''(\overline{\eta})+\zeta(\overline{\eta})\zeta'''(\overline{\eta})]}{r^2(\zeta'(\overline{\eta})-1)}-\frac{\zeta(\overline{\eta})\zeta''(\overline{\eta})[\zeta'(\overline{\eta})-1+\zeta"(\overline{\eta})\zeta''(\overline{\eta})]}{r^2(\zeta'(\overline{\eta})-1)^2}\Bigg\}\frac{d^2\overline{\eta}}{dr^2}&
  \end{eqnarray*}
converges as $r\to0$ by applying \eqref{e5.14} again. Next, using the identity
  $$
   \overline{\eta}(r)=1+\frac{\overline{\eta}''(0)}{2}r^2+\int^r_0\Bigg[\frac{1}{6}(r^3-s^3)+\frac{1}{2}s^2(r-s)-\frac{s}{2}(r^2-s^2)\Bigg]\overline{\eta}^{(4)}(s)ds
  $$
and a similar computation given as above, we conclude that $v\in C^4([0,1])$ and $v^{(4)}(0)=0$. The proof was done. $\Box$\\

\vspace{40pt}

\section{Schauder's fix point theorem of \eqref{e4.7} for negative pair}

 In this section, we will use Schauder's fix point theorem to prove the following local existence result for some negative constant $\lambda'''$.

 \begin{theo}\label{t6.1}
   Letting $2\leq n\leq5$ and $\theta>0$, there exists a negative constant $\lambda'''$ such that \eqref{e4.7} admits a local positive convex solution $\zeta\in C^{3,1}([1,r_0])\cap C^\infty((1,r_0])$ for some $r_0>1$. Furthermore, the solution satisfies
     \begin{equation}\label{e6.1}
      \zeta(1)=0, \ \ \zeta'(1)=2, \ \ \zeta''(1)=\frac{4(n+2)^2\theta+(2n^2-24n+104)}{n^2-2n+24}.
     \end{equation}
 \end{theo}

  \noindent\textbf{Proof.} At first, defining a family of functions
    \begin{eqnarray*}
     &\Gamma_{\eta_0, \alpha, \beta, \gamma, \sigma}\equiv\Big\{\varphi\in C^3([1,\eta_0])\Big|\ \varphi(1)=0, \varphi'(1)=2,\ \ \varphi''(1)=\alpha,\ \ \varphi'''(1)=\beta&\\
     &\varphi(\eta)\in[0,1],\ \ \varphi'(\eta)\in[2-\sigma,2+\sigma], \ \  \varphi''(\eta)\in[\alpha-\sigma, \alpha+\sigma],\ \ \forall \eta\in[1,\eta_0]&\\
     &\frac{\varphi'''(\eta)-\varphi'''(1)}{\eta-1}\in[\gamma-1,\gamma+1], \ \ \varphi'''(\eta)\in[\beta-\sigma,\beta+\sigma], \ \ \forall \eta\in(1,\eta_0]\Big\}&
    \end{eqnarray*}
  for any $\alpha, \beta, \gamma-1<\gamma+1\in{\mathbb{R}}, \eta_0>1, \sigma>0$, it's clear that closed convex subset of $C^3([1,\eta_0])$ endowed with norm $||\cdot||_{C^3([1,\eta_0])}$. Given any $\varphi\in\Gamma_{\eta_0, \alpha, \beta, \gamma, \sigma}$, let's introduce a mapping $T\varphi\equiv\zeta$ by
    \begin{equation}\label{e6.2}
     \begin{cases}
       \zeta'=(\theta+1)\eta^{-1}\varphi+n\eta(\eta-1)\Big\{\big[n\theta-(n-1)\big]\eta-\big[n\theta-1\big]\Big\}\varphi^{-1}\\
       \ \ \ \  +\Big\{\big[2n\theta-(2n-1)\big]\eta-\big[2n\theta-1\big]\Big\}+\lambda(\varphi,\eta_0)\frac{\eta^2}{\varphi}\exp\int^\eta_{\eta_0}\frac{s+1}{\varphi(s)}ds, &\forall \eta\in(1,\eta_0]\\
       \zeta(1)=0,
     \end{cases}
    \end{equation}
 where $\lambda(\varphi,\eta_0)$ is a constant being chosen such that
   \begin{equation}\label{e6.3}
     \lim_{\eta\to1^+}\lambda(\varphi,\eta_0)\frac{1}{\varphi(\eta)}\exp\int^\eta_{\eta_0}\frac{s+1}{\varphi(s)}ds=4+\frac{n(n-2)}{2}.
   \end{equation}
 It's clear that the solution $\zeta\in C^4((1,\eta_0])$. Next, we want to show that for appropriate chosen $\alpha, \beta, \gamma\in{\mathbb{R}}$ and $\sigma,\eta_0$ small, $T$ is a continuous and compact mapping from convex set $\Gamma_{\eta_0, \alpha, \beta, \gamma, \sigma}$ to itself. Then by Schauder's fix point theorem \cite{C}, there exists a fix point belonging to $\Gamma_{\eta_0, \alpha, \beta, \gamma, \sigma}$ and satisfying \eqref{e4.7}. Our arguments are divided into several crucial lemmas.

 \begin{lemm}\label{l6.1}
   For any $\varphi\in\Gamma_{\eta_0, \alpha, \beta, \gamma, \sigma}$, we have
     \begin{equation}\label{e6.4}
       \Phi(\eta)\equiv\frac{1}{\varphi(\eta)}\exp\int^\eta_{\eta_0}\frac{s+1}{\varphi(s)}ds\leq\frac{\alpha+\sigma}{8}\frac{\eta_0-1+\frac{4}{\alpha+\sigma}}{\eta_0-1},\ \ \forall \eta\in(1,\eta_0).
     \end{equation}
   Moreover, if $\alpha-\sigma\geq1$, then $\Phi(\cdot)$ is a monotone non-increasing function in $\eta$ and hence
     \begin{equation}\label{e6.5}
       \Phi(\eta)\leq\lim_{\eta\to1^+}\Phi(\eta)=\frac{4}{\lambda(\varphi,\eta_0)},\ \ \forall \eta\in[1,\eta_0].
     \end{equation}
 \end{lemm}

\noindent\textbf{Proof.} Since $\alpha-\sigma\leq\varphi''\leq\alpha+\sigma$ for $\varphi\in\Gamma_{\eta_0, \alpha, \beta, \gamma, \sigma}$, it's inferred from the initial conditions of $\varphi$ that
   \begin{equation}\label{e6.6}
    2(s-1)+\frac{\alpha-\sigma}{2}(s-1)^2\leq\varphi(s)\leq 2(s-1)+\frac{\alpha+\sigma}{2}(s-1)^2, \ \ \forall s\in(1,\eta_0).
   \end{equation}
Therefore, \eqref{e6.4} follows from the integration
   $$
    \int^\eta_{\eta_0}\frac{2}{2(s-1)+\frac{\alpha+\sigma}{2}(s-1)^2}ds=\log\Bigg(\frac{\eta-1}{\eta-1+\frac{4}{\alpha+\sigma}}\Bigg)-\log\Bigg(\frac{\eta_0-1}{\eta_0-1+\frac{4}{\alpha+\sigma}}\Bigg).
   $$
To show the monotonicity of $\Phi(\cdot)$, we need only calculate the derivative of $\Phi$ by
  \begin{eqnarray*}
   \Phi'(\eta)&=&\frac{-\varphi'+(\eta+1)}{\varphi^2}\exp\int^\eta_{\eta_0}\frac{s+1}{\varphi(s)}ds\\
    &\leq&\frac{\eta-1}{\varphi^2}(1-\beta/2)\leq0, \ \ \forall \eta\in[1,\eta_0].
  \end{eqnarray*}
The proof was done. $\Box$\\

\begin{lemm}\label{l6.2}
  For any $\varphi\in\Gamma_{\eta_0, \alpha, \beta, \gamma, \sigma}$, we have
    \begin{equation}\label{e6.7}
       \frac{32+4n(n-2)}{\alpha+\sigma}\frac{\eta_0-1}{\eta_0-1+\frac{4}{\alpha+\sigma}}\leq\lambda(\varphi,\eta_0)\leq(\eta_0-1)\Bigg(4+\frac{n(n-2)}{2}\Bigg)\Bigg(2+\frac{\alpha+\sigma}{2}(\eta_0-1)\Bigg)e^{\frac{\eta_0-1}{2}}.
    \end{equation}
\end{lemm}

\noindent\textbf{Proof.} Noting that for $\varphi\in\Gamma_{\eta_0, \alpha, \beta, \gamma, \sigma}$, there holds
  \begin{eqnarray}\nonumber\label{e6.8}
   \frac{1}{\varphi}\exp\int^\eta_{\eta_0}\frac{s+1}{\varphi(s)}ds&\geq&\frac{1}{2(\eta-1)+\frac{\alpha+\sigma}{2}(\eta-1)^2}\exp\int^\eta_{\eta_0}\frac{s+1}{2(s-1)}ds\\
    &\geq&\frac{1}{(\eta_0-1)\Big(2+\frac{\alpha+\sigma}{2}(\eta_0-1)\Big)}e^{\frac{1-\eta_0}{2}}.
  \end{eqnarray}
So, \eqref{e6.7} follows from \eqref{e6.4} and \eqref{e6.8}. $\Box$\\

\begin{lemm}\label{l6.3}
 $\lambda(\cdot,\eta_0)$ is a continuous function in $C^3([1,\eta_0])$ for each fixed $\eta_0$.\\
\end{lemm}

\noindent\textbf{Proof.} Using
   $$
    \varphi(s)=2(s-1)+\int^s_1(s-\tau)\varphi''(\tau)d\tau,
   $$
one has
   \begin{eqnarray*}
    F(\varphi)&\equiv&\frac{1}{\varphi(\eta)}\exp\int^\eta_{\eta_0}\frac{s+1}{\varphi(s)}ds\\
     &=&\frac{G(\varphi,\eta)}{\eta-1}\exp\int^\eta_{\eta_0}\frac{s+1}{s-1}G(\varphi,s)ds,
   \end{eqnarray*}
where
    $$
     G(\varphi,s)\equiv\frac{1}{2+\frac{1}{s-1}\int^s_1(s-\tau)\varphi''(\tau)d\tau}\in C(C^3([1,\eta_0])\times[1,\eta_0])
    $$
Next, we show that for $\varphi_\varepsilon$ tends to $\varphi$ in $C^3([1,\eta_0])$ as $\varepsilon\to0$, there holds
  $$
    H(\varphi_\varepsilon,\eta)\equiv\frac{1}{\eta-1}\exp\int^\eta_{\eta_0}\frac{s+1}{s-1}G(\varphi_\varepsilon,s)ds
  $$
tends to $H(\varphi,\eta)$ uniformly in $\eta\in[1,\eta_0]$ as $\varepsilon\to0$. In fact,
  \begin{eqnarray*}
    \frac{H(\varphi_\varepsilon,\eta)}{H(\varphi,\eta)}&=&\exp\int^\eta_{\eta_0}\frac{s+1}{s-1}(G(\varphi_\varepsilon,s)-G(\varphi,s))ds\\
      &\leq&\exp\Bigg\{o(\varepsilon)\int^\eta_{\eta_0}\frac{s+1}{s-1}ds\Bigg\}=\exp\Bigg\{o(\varepsilon)\Bigg[(\eta-\eta_0)+2\log\Bigg(\frac{\eta-1}{\eta_0-1}\Bigg)\Bigg\}\\
      &=&\Bigg(\frac{\eta-1}{\eta_0-1}\Bigg)^{2o(\varepsilon)}\exp\Big\{o(\varepsilon)(\eta-\eta_0)\Big\}=1+o(\varepsilon),
  \end{eqnarray*}
 where $o(\varepsilon)$ is a small quantity depending only on $\varepsilon$ as $\varepsilon\to0$. Similarly
    $$
     \frac{H(\varphi,\eta)}{H(\varphi_\varepsilon,\eta)}\leq 1+o(\varepsilon).
    $$
Another hand, noting that
   \begin{eqnarray*}
     H(\varphi,\eta)&\geq&\frac{1}{\eta-1}\exp\int^\eta_{\eta_0}\frac{s+1}{2(s-1)}ds\\
      &=&\frac{1}{\eta-1}\frac{\eta-1}{\eta_0-1}e^{\frac{\eta-\eta_0}{2}}=\frac{1}{\eta_0-1}e^{\frac{\eta-\eta_0}{2}},
   \end{eqnarray*}
the uniformly convergence of $H(\varphi_\varepsilon,s)$ to $H(\varphi,s)$ follows from
   $$
    |H(\varphi_\varepsilon,\eta)-H(\varphi,\eta)|=\frac{\Big|\frac{H(\varphi_\varepsilon,\eta)}{H(\varphi,\eta)}-1\Big|}{H(\varphi,\eta)}.
   $$
$\Box$\\

\begin{lemm}\label{l6.4}
  For any $\varphi\in\Gamma_{\eta_0, \alpha, \beta, \gamma, \sigma}, \alpha-\sigma\geq1$, there holds
     \begin{equation}\label{e6.9}
       0\leq\zeta(\eta)\leq1, \ \ \forall[1,\eta_0]
     \end{equation}
for solution $\zeta=T\varphi$, provided $\zeta_0$ is small.
\end{lemm}

\noindent\textbf{Proof.} Using \eqref{e6.2} and Lemma \ref{l6.1}, one gets that
  \begin{eqnarray*}
    \zeta'(\eta)&\leq&(\theta+1)+\max\Big\{-2,\big[(2n\theta-(2n-1))\eta_0-(2n\theta-1)\big]\Big\}\\
     &&+\frac{n\eta_0}{2}\Big\{\big[n\theta-(n\theta-1)\big]\eta_0-\big[n\theta-1\big]\Big\}+4\eta_0^2, \ \ \forall \eta\in(1,\eta_0).
  \end{eqnarray*}
Integrating over $\eta$, we obtain that
  \begin{eqnarray*}
   \zeta(\eta)&\leq&C_{\theta}(\eta-1)\leq1, \ \ \forall \eta\leq\eta_0
  \end{eqnarray*}
by choosing $\eta_0$ small. $\Box$\\

\begin{lemm}\label{l6.5}
 For any $\varphi\in\Gamma_{\eta_0, \alpha, \beta, \gamma, \sigma}$, the solution $\zeta\equiv T\varphi$ satisfies that
   \begin{equation}\label{e6.10}
     \lim_{\eta\to1^+}\zeta'(\eta)=2,\ \ \zeta'(\eta)\in[2-\sigma,2+\sigma], \ \ \forall \eta\in[1,\eta_0],
   \end{equation}
 provided $\eta_0-1$ is chosen small with respect to $\sigma$. As a result, $\zeta\in C^1([1,\eta_0])$ and $\zeta'(1)=2$.
\end{lemm}

\noindent\textbf{Proof.} By \eqref{e6.2} and
   $$
     2(\eta-1)+\frac{\alpha-\sigma}{2}(\eta-1)^2\leq\varphi\leq 2(\eta-1)+\frac{\alpha+\sigma}{2}(\eta-1)^2, \ \ \forall \eta\in[1,\eta_0],
   $$
it follows from \eqref{e6.2} and \eqref{e6.3} that
   $$
    \lim_{\eta\to1^+}\zeta'(\eta)=2
   $$
and
   $$
    \zeta'(\eta)=2+o(\eta_0-1)\in[2-\sigma,2+\sigma], \ \ \forall \eta\in(1,\eta_0].
   $$
Therefore, by Lagrange's intermediate value theorem,
   $$
    \zeta'(1)=\lim_{\eta\to1^+}\frac{\zeta(\eta)-\zeta(1)}{\eta-1}=\lim_{\xi\to1^+}\zeta'(\xi), \ \ \xi\in(1,\eta)
   $$
exists and equals to $2$. And thus $\zeta\in C^1([1,\eta_0])$. $\Box$\\

\begin{lemm}\label{l6.6}
 For
   $$
     \alpha\equiv\frac{4(n+2)^2\theta+(2n^2-24n+104)}{n^2-2n+24}
   $$
 and $\varphi\in\Gamma_{\eta_0, \alpha, \beta, \gamma, \sigma}$, the solution $\zeta\equiv T\varphi$ satisfies that
   \begin{equation}\label{e6.11}
     \lim_{\eta\to1^+}\zeta''(\eta)=\alpha,\ \ \zeta''(\eta)\in[\alpha-\sigma,\alpha+\sigma], \ \ \forall \eta\in[1,\eta_0],
   \end{equation}
 provided $\eta_0-1$ is chosen smaller with respect to $\sigma$. As a result, $\zeta\in C^2([1,\eta_0])$ and $\zeta''(1)=\alpha$.
\end{lemm}

\noindent\textbf{Proof.} Imposing again
   \begin{eqnarray*}
    &2(s-1)+\frac{\alpha-\sigma}{2}(s-1)^2\leq\varphi(s)\leq 2(s-1)+\frac{\alpha+\sigma}{2}(s-1)^2,&\\
    &2+(\alpha-\sigma)(s-1)\leq\varphi'(s)\leq2+2(\alpha+\sigma)(s-1),&\\
    &\alpha-\sigma\leq\varphi''(s)\leq\alpha+\sigma, \ \ \forall s\in[1,\eta_0]&
   \end{eqnarray*}
together with the computation
  \begin{eqnarray}\nonumber\label{e6.12}
    &\zeta''=-(\theta+1)\eta^{-2}\varphi+(\theta+1)\eta^{-1}\varphi'+\Big[2n\theta-(2n-1)\Big]&\\ &+n\Big\{3\big[n\theta-(n-1)\big]\eta^2-2n(2\theta-1)\eta+(n\theta-1)\Big\}\varphi^{-1}&\\ \nonumber
     &-n\eta(\eta-1)\Big\{\big[n\theta-(n-1)\big]\eta-(n\theta-1)\Big\}\varphi^{-2}\varphi'+\frac{\lambda\eta}{\varphi}\exp\int^\eta_{\eta_0}\frac{s+1}{\varphi(s)}ds\Big[\frac{2\varphi-\eta\varphi'+\eta(\eta+1)}{\varphi}\Big],&
  \end{eqnarray}
one gets
    \begin{eqnarray*}
     \zeta''(\eta)&=&2(n+1)\theta-(2n-3)+J_1+J_2\frac{\lambda\eta}{\varphi}\exp\int^\eta_{\eta_0}\frac{s+1}{\varphi(s)}ds+o(\eta_0-1)\\
       &=&2(n+1)\theta-(2n-3)+\frac{n[n\theta-(2n-3)]}{2}+\frac{n(n-2)}{8}\alpha\\
       &&+\frac{[8+n(n-2)](5-\alpha)}{4}+o(\eta_0-1)=\alpha+o(\eta_0-1),\ \ \forall \eta\in[1,\eta_0],
    \end{eqnarray*}
 where
  \begin{eqnarray*}
    J_1&\equiv&\frac{n\varphi\Big\{3\big[n\theta-(n-1)\big]\eta^2-2n(2\theta-1)\eta+(n\theta-1)\Big\}}{\varphi^2}\\
      &&-\frac{n\eta(\eta-1)\Big\{[n\theta-(n-1)]\eta-(n\theta-1)\Big\}\varphi'}{\varphi^2}\\
      &=&\frac{n\Big\{6[n\theta-(n-1)]\xi-2n(2\theta-1)\Big\}\varphi-n\xi(\xi-1)\Big\{[n\theta-(n-1)]\xi-(n\theta-1)\Big\}\varphi''}{2\varphi\varphi'}\\
      &=&\frac{n[n\theta-(2n-3)]}{2}+\frac{n(n-2)}{8}\alpha+o(\eta_0-1), \ \ \xi\in(1,\eta)
  \end{eqnarray*}
and
   \begin{eqnarray*}
     J_2&\equiv&\frac{2\varphi-\eta\varphi'+\eta(\eta+1)}{\varphi}\\
       &=&\frac{\varphi'-\xi\varphi''+2\xi+1}{\varphi'}=\frac{5-\alpha}{2}+o(\eta_0-1), \ \ \xi\in(1,\eta)
   \end{eqnarray*}
has been used by Cauchy's intermediate value theorem. So the lemma can be proven as above without difficulty. $\Box$\\

\begin{lemm}\label{l6.7}
 For
 \begin{eqnarray*}
   \alpha&=&\alpha_0(n,\theta)\equiv\frac{4(n+2)^2\theta+(2n^2-24n+104)}{n^2-2n+24},\\ \beta&=&\beta_0(n,\theta)\equiv\frac{48(n+2)(n-2)\theta+6(n-2)(9n-8)+528}{96+2n(n-2)}\\
    &&+\frac{6[n(n-2)+12]\alpha_0^2-3[4(n+2)(n-2)\theta+(n-2)(13n-4)+144]\alpha_0}{96+2n(n-2)}
 \end{eqnarray*}
and  $\varphi\in\Gamma_{\eta_0, \alpha, \beta, \gamma, \sigma}$, the solution $\zeta\equiv T\varphi$ satisfies that
   \begin{equation}\label{e6.13}
     \lim_{\eta\to1^+}\zeta'''(\eta)=\beta, \ \ \zeta'''(\eta)\in[\beta-\sigma,\beta+\sigma],\ \  \forall \eta\in[1,\eta_0],
   \end{equation}
 provided $\eta_0-1$ is chosen smaller with respect to $\sigma$. As a result, $\zeta\in C^3([1,\eta_0])$ and $\zeta'''(1)=\beta$.
\end{lemm}

\noindent\textbf{Proof.} In fact, using
  \begin{eqnarray}\nonumber\label{e6.14}
    &\varphi(s)\geq2(s-1)+\frac{\alpha}{2}(s-1)^2+\frac{\beta}{6}(s-1)^3+\frac{\gamma-1}{240}(s-1)^4&\\ \nonumber
    &\varphi(s)\leq2(s-1)+\frac{\alpha}{2}(s-1)^2+\frac{\beta}{6}(s-1)^3+\frac{\gamma+1}{240}(s-1)^4&\\
    \nonumber
    &\varphi'(s)\geq2+\alpha(s-1)+\frac{\beta}{2}(s-1)^2+\frac{\gamma-1}{60}(s-1)^3&\\
    &\varphi'(s)\leq2+\alpha(s-1)+\frac{\beta}{2}(s-1)^2+\frac{\gamma+1}{60}(s-1)^3&\\ \nonumber
    &\alpha+\beta(s-1)+\frac{\gamma-1}{20}(s-1)^2\leq\varphi''(s)\leq \alpha+\beta(s-1)+\frac{\gamma+1}{20}(s-1)^2&\\ \nonumber
    &\beta+\frac{\gamma-1}{10}(s-1)\leq\varphi'''(s)\leq\beta+\frac{\gamma+1}{10}(s-1),&
  \end{eqnarray}
we get the formula for third derivative
  \begin{eqnarray*}
    \zeta'''&=&K_1+K_2+RJ_2+\eta RR'\\
    &=&\frac{2[n(n-2)+12]\alpha^2-[4(n+2)(n-2)\theta+(n-2)(13n-4)+144]\alpha}{16}\\
    &&-\Bigg[1+\frac{n(n-2)}{24}\Bigg]\beta+(n+2)(n-2)\theta+\frac{(n-2)(9n-8)}{8}+11+o(\eta_0-1),
  \end{eqnarray*}
where
   \begin{eqnarray*}
     K_1&\equiv&2(\theta+1)\eta^{-3}\varphi-2(\theta+1)\eta^{-2}\varphi'+(\theta+1)\eta^{-1}\varphi''\\
     &=&-4(\theta+1)+(\theta+1)\alpha+o(\eta_0-1),\\
     R&\equiv&\frac{\lambda}{\varphi}\exp\int^\eta_{\eta_0}\frac{s+1}{\varphi(s)}ds=4+\frac{n(n-2)}{2}+o(\eta_0-1)
   \end{eqnarray*}
together with
    \begin{eqnarray*}
      K_2&\equiv&\frac{\partial J_1}{\partial\eta}=\frac{n\varphi^2\Big\{6[n\theta-(n-1)]\eta-2n(2\theta-1)\Big\}}{\varphi^3}\\
      &&-\frac{2n\varphi\varphi'\Big\{3[n\theta-(n-1)]\eta^2-2n(2\theta-1)\eta+(n\theta-1)\Big\}}{\varphi^3}\\
      &&-\frac{n\eta(\eta-1)\Big\{[n\theta-(n-1)]\eta-(n\theta-1)\Big\}\Big[\varphi\varphi''-2(\varphi')^2\Big]}{\varphi^3}\\
      &=&n[n\theta-(n-1)]+\frac{n(n-2)}{12}\beta+\frac{3n\xi(\xi-1)\Big\{[n\theta-(n-1)]\xi-(n\theta-1)\Big\}\varphi'\varphi''}{3\varphi^3\varphi'}\\
      &&-\frac{3n\Big\{3[n\theta-(n-1)]\xi^2-2n(2\theta-1)\xi+(n\theta-1)\Big\}\varphi\varphi''}{3\varphi^2\varphi'}+o(\eta_0-1)\\
      &=&n[n\theta-(n-1)]+\frac{n(n-2)}{12}\beta-\frac{n(n-2)\alpha^2+4n[n\theta-(2n-3)]\alpha}{16}+o(\eta_0-1)
    \end{eqnarray*}
and
  \begin{eqnarray*}
    R'&\equiv&\Bigg(-\frac{\varphi'}{\varphi}+\frac{\eta+1}{\varphi}\Bigg)J_2+\frac{\partial J_2}{\partial\eta}\\
     &=&\frac{-3\varphi\varphi'+2\eta(\varphi')^2-3\eta(\eta+1)\varphi'+(4\eta+3)\varphi+\eta(\eta+1)^2-\eta\varphi\varphi''}{\varphi^2}\\
     &=&\frac{-4\varphi\varphi''-(\varphi')^2+3\xi\varphi'\varphi''-2\xi\varphi'-3\xi(\xi+1)\varphi''+4\varphi+(3\xi+1)(\xi+1)-\xi\varphi\varphi'''}{2\varphi\varphi'}\\
     &=&-\alpha+1-\frac{\beta}{4}+\frac{-2\varphi'\varphi''+3\varphi'\varphi''+3\chi(\varphi'')^2+3\chi\varphi'\varphi'''-2\varphi'-2\chi\varphi''}{2(\varphi')^2+2\varphi\varphi''}\\
     &&+\frac{-(6\chi+3)\varphi''-3\chi(\chi+1)\varphi'''+6\chi+4}{2(\varphi')^2+2\varphi\varphi''}\\
     &=&-\frac{\beta}{4}+\frac{3\alpha^2-17\alpha+14}{8}+o(\eta_0-1)
  \end{eqnarray*}
have been used by Cauchy's intermediate value theorem for $\xi\in(1,\eta)$ and $\chi\in(1,\xi)$. Noting that for $\beta=\beta_0(n,\theta)$, we have
   \begin{equation}\label{e6.15}
     \zeta'''(\eta)=\beta_0(n,\theta)+o(\eta_0-1), \ \ \forall\eta\in(1,\eta_0]
   \end{equation}
and complete the proof of \eqref{e6.13} and $\zeta\in C^3([1,\eta_0])$. $\Box$\\

\begin{lemm}\label{l6.8}
 Suppose that $2\leq n\leq5$. For
  \begin{equation}\label{e6.16}
      \alpha=\alpha_0(n,\theta),\ \ \beta=\beta_0(n,\theta), \ \ \gamma=\frac{48\gamma_{\alpha,\beta}}{80+n(n-2)}
  \end{equation}
and  $\varphi\in\Gamma_{\eta_0, \alpha, \beta, \gamma, \sigma}$, the solution $\zeta\equiv T\varphi$ satisfies that
   \begin{equation}\label{e6.17}
     \zeta^{(4)}(\eta)\in[\gamma-1,\gamma+1],\ \  \forall \eta\in(1,\eta_0],
   \end{equation}
 provided $\sigma>0$ is chosen small and then $\eta_0-1$ is also chosen smaller with respect to $\sigma$, where
    \begin{eqnarray*}
     \gamma_{\alpha,\beta}&\equiv&-\frac{17[8+n(n-2)]}{96}\alpha^3+\frac{9n[n\theta-(2n-3)]+53[8+n(n-2)]}{48}\alpha^2\\
      &&-\frac{216n[n\theta-(n-1)]+1021[8+n(n-2)]}{288}\alpha+\frac{37[8+n(n-2)]}{16}\\
      &&-\frac{12n[n\theta-(2n-3)]+61[8+n(n-2)]}{48}\beta+\frac{112-n(n-2)}{48}\alpha\beta.
    \end{eqnarray*}
  As a result, $\zeta$ lies in a bounded set of H\"{o}lder space $C^{3,1}([1,\eta_0])$.
\end{lemm}

\noindent\textbf{Proof.} Taking once more derivative on $\zeta$ yields that
  \begin{equation}\label{e6.18}
    \zeta^{(4)}=P+Q+2RR'+\eta RR'\frac{-\varphi'+\eta+1}{\varphi}+\eta R\frac{\partial R'}{\partial\eta},
  \end{equation}
where
   \begin{eqnarray*}
    P&\equiv&\frac{\partial K_1}{\partial\eta}=-6(\theta+1)\eta^{-4}\varphi+6(\theta+1)\eta^{-3}\varphi'-3(\theta+1)\eta^{-2}\varphi''\\
    &&+(\theta+1)\eta^{-1}\varphi'''=-12(\theta+1)-3(\theta+1)\alpha+(\theta+1)\beta+o(\eta_0-1),
   \end{eqnarray*}
  \begin{eqnarray*}
    Q&\equiv&\frac{\partial K_2}{\partial\eta}=6n[n\theta-(n-1)]\varphi^{-1}-3n\Big\{6[n\theta-(n-1)]\eta-2n(2\theta-1)\Big\}\frac{\varphi'}{\varphi^2}\\
    &&-3n\Big\{3[n\theta-(n-1)]\eta^2-2n(2\theta-1)\eta+(n\theta-1)\Big\}\Bigg[\frac{\varphi''}{\varphi^2}-2\frac{(\varphi')^2}{\varphi^3}\Bigg]\\
    &&-n\eta(\eta-1)\Big\{[n\theta-(n-1)\Big]\eta-(n\theta-1)\Big\}\Bigg[\frac{\varphi'''}{\varphi^2}-6\frac{\varphi'\varphi''}{\varphi^3}+6\frac{(\varphi')^3}{\varphi^4}\Bigg]\\
    &=&\frac{n(n-2)}{16}\frac{\varphi'''(\chi)-\varphi'''(1)}{\chi-1}+\frac{3n[n\theta-(2n-3)]\alpha^2-12n[n\theta-(n-1)]\alpha}{16}\\
    &&-\frac{4n[n\theta-(2n-3)]\beta+5n(n-2)\alpha\beta}{16}+o(\eta_0-1)
   \end{eqnarray*}
and
   \begin{eqnarray*}
    \frac{\partial R'}{\partial\eta}&=&\frac{5\varphi(\varphi')^2-4\varphi^2\varphi''+5\eta\varphi\varphi'\varphi''-(10\eta+6)\varphi\varphi'-3\eta(\eta+1)\varphi\varphi''+4\varphi^2}{\varphi^3}\\
     &&+\frac{(\eta+1)(3\eta+1)\varphi-\eta\varphi^2\varphi'''-4\eta(\varphi')^3+6\eta(\eta+1)(\varphi')^2-2\eta(\eta+1)^2\varphi'}{\varphi^3}\\
     &=&-\frac{\varphi'''(\chi)-\varphi'''(1)}{6(\chi-1)}-\frac{23}{12}\beta+\frac{61}{48}\alpha^2+\frac{11}{24}\alpha\beta-\frac{9}{48}\alpha^3-\frac{65}{72}\alpha+\frac{1}{4}+o(\eta_0-1)
   \end{eqnarray*}
by \eqref{e6.14} and Cauchy's intermediate value theorem, where $\chi\in(1,\eta)$. Therefore, one obtain that
  \begin{equation}\label{e6.19}
   \zeta^{(4)}(\eta)=\gamma_{\alpha,\beta}-\Bigg[\frac{2}{3}+\frac{n(n-2)}{48}\Bigg]\frac{\varphi'''(\chi)-\varphi'''(1)}{\chi-1}\in[\gamma-1,\gamma+1], \ \ \forall\eta\in[1,\eta_0]
  \end{equation}
upon \eqref{e6.16} and choosing $\sigma, \eta_0$ small. The proof was done. $\Box$\\

Summarizing the above lemmas shows that $T$ is a continuous and compact mapping from $\Gamma_{\eta_0, \alpha, \beta, \gamma, \sigma}$ to $\Gamma_{\eta_0, \alpha, \beta, \gamma, \sigma}$.

\begin{prop}\label{p6.1}
  Under the assumption of Lemma \ref{l6.1}-\ref{l6.8}, the mapping $T: \Gamma_{\eta_0, \alpha, \beta, \gamma, \sigma}\to\Gamma_{\eta_0, \alpha, \beta, \gamma, \sigma}$ is continuous and compact.
\end{prop}

\noindent\textbf{Proof.} By continuity of $F(\cdot)$ and $\lambda(\cdot,\eta_0)$ on $C^3([1,\eta_0])$ proven in Lemma \ref{l6.3}, $T$ is a continuous mapping from $C^3([1,\eta_0])$ to $C^1([1,\eta_0])$. To show the continuity of $T$ from $\Gamma_{\eta_0, \alpha, \beta, \gamma, \sigma}$ to $C^3([1,\eta_0])$, we use the formula
  $$
   \zeta'''=K_1(\varphi)+K_2(\varphi)+R(\varphi)J_2(\varphi)+\eta R(\varphi)R'(\varphi)
  $$
of third derivative of $\zeta$. Taking a sequence $\varphi_\varepsilon$ approaching to $\varphi$ in $\Gamma_{\eta_0, \alpha, \beta, \gamma, \sigma}$, it's clear that
  \begin{equation}\label{e6.20}
    K_1(\varphi_\varepsilon)\to K_1(\varphi), \ \ \mbox{as } \varepsilon\to0
  \end{equation}
and
   \begin{equation}\label{e6.21}
     R(\varphi_\varepsilon)\to R(\varphi), \ \ \mbox{as } \varepsilon\to0
   \end{equation}
by Lemma \ref{l6.3}. Secondly, we have
   \begin{equation}\label{e6.22}
     J_2(\varphi_\varepsilon)\to J_2(\varphi), \ \ \mbox{as } \varepsilon\to0
   \end{equation}
by
   \begin{eqnarray}\nonumber\label{e6.23}
    &\Big|\frac{-\eta\varphi_\varepsilon'+\eta(\eta+1)}{\varphi_\varepsilon}-\frac{-\eta\varphi'+\eta(\eta+1)}{\varphi}\Big|=\frac{\big|-\eta\varphi\varphi_\varepsilon'+\eta(\eta+1)\varphi+\eta\varphi'\varphi_\varepsilon-\eta(\eta+1)\varphi_\varepsilon\big|}{\varphi_\varepsilon\varphi}&\\ \nonumber
    &\leq\frac{\eta-1}{\varphi_\varepsilon\varphi}\Big|-\varphi\varphi_\varepsilon'+(2\eta+1)\varphi+\eta(\eta+1)\varphi'+\varphi'\varphi_\varepsilon-(2\eta+1)\varphi_\varepsilon-\eta(\eta+1)\varphi_\varepsilon'\Big|&\\ \nonumber
    &+\frac{\eta(\eta-1)}{\varphi_\varepsilon}|\varphi_\varepsilon''-\varphi''|+\frac{\eta(\eta-1)}{\varphi_\varepsilon\varphi}|\varphi_\varepsilon-\varphi||\varphi''|&\\ \nonumber
    &\leq\frac{\eta-1}{\varphi_\varepsilon}|\varphi_\varepsilon'-\varphi'|+\frac{\eta-1}{\varphi_\varepsilon\varphi}|\varphi_\varepsilon-\varphi||\varphi'|+\frac{(2\eta+1)(\eta-1)}{\varphi_\varepsilon\varphi}|\varphi_\varepsilon'-\varphi'|+\frac{\eta(\eta+1)(\eta-1)}{\varphi_\varepsilon\varphi}|\varphi_\varepsilon'-\varphi|&\\ \nonumber
    &+\frac{\eta(\eta-1)}{\varphi_\varepsilon}|\varphi_\varepsilon''-\varphi''|+\frac{\eta(\eta-1)}{\varphi_\varepsilon\varphi}|\varphi_\varepsilon-\varphi||\varphi''|&\\ \nonumber
    &\leq\frac{1}{2+\frac{\beta}{4}(\eta-1)}o(\varepsilon)+\frac{1}{(2+\frac{\beta}{4}(\eta-1))^2}|\varphi'|o(\varepsilon)+\frac{2\eta+1}{(2+\frac{\beta}{4}(\eta-1))^2}o(\varepsilon)+\frac{\eta(\eta+1)}{(2+\frac{\beta}{4}(\eta-1))^2}o(\varepsilon)&\\
    &+\frac{\eta}{2+\frac{\beta}{4}(\eta-1)}o(\varepsilon)+\frac{\eta}{(2+\frac{\beta}{4}(\eta-1))^2}|\varphi''|o(\varepsilon).&
   \end{eqnarray}
Noting that
   \begin{eqnarray}\label{e6.24}
     \Bigg|\frac{-\varphi_\varepsilon'+\eta+1}{\varphi_\varepsilon}-\frac{-\varphi'+\eta+1}{\varphi}\Bigg|=o(\varepsilon)
   \end{eqnarray}
by \eqref{e6.23}, and
   \begin{eqnarray}\nonumber\label{e6.25}
     \frac{\partial J_2(\varphi_\varepsilon)}{\partial\eta}&=&\frac{-3\varphi_\varepsilon\varphi_\varepsilon'-\eta\varphi_\varepsilon\varphi_\varepsilon''+(2\eta+1)\varphi_\varepsilon+\eta(\varphi_\varepsilon')^2-\eta(\eta+1)\varphi_\varepsilon'}{\varphi_\varepsilon^2}\\ &=&\frac{\partial J_2(\varphi)}{\partial\eta}+o(\varepsilon)
   \end{eqnarray}
by
    \begin{eqnarray}\nonumber\label{e6.26}
      &\displaystyle \Bigg|\frac{1}{\varphi_\varepsilon^2}-\frac{1}{\varphi^2}\Bigg|=\Bigg|\frac{(\varphi_\varepsilon+\varphi)(\varphi_\varepsilon-\varphi)}{\varphi_\varepsilon^2\varphi^2}\Bigg|\leq C\Bigg|\frac{\varphi_\varepsilon-\varphi}{\varphi_\varepsilon\varphi^2}\Bigg|=C\Bigg|\frac{\varphi_\varepsilon'-\varphi'}{\varphi_\varepsilon'\varphi^2+2\varphi_\varepsilon\varphi\varphi'}\Bigg|&\\
       &\displaystyle =C\Bigg|\frac{\varphi_\varepsilon''-\varphi''}{\varphi_\varepsilon''\varphi^2+4\varphi_\varepsilon'\varphi\varphi'+2\varphi_\varepsilon(\varphi')^2+2\varphi_\varepsilon\varphi\varphi''}\Bigg|&\\ \nonumber
       &\displaystyle =C\Bigg|\frac{\varphi_\varepsilon'''-\varphi'''}{6\varphi_\varepsilon'(\varphi')^2+6\varphi_\varepsilon\varphi'\varphi''+6\varphi_\varepsilon'\varphi\varphi''+6\varphi_\varepsilon''\varphi\varphi'+2\varphi_\varepsilon\varphi'''+\varphi_\varepsilon'''\varphi^2}\Bigg|=o(\varepsilon),&
    \end{eqnarray}
we conclude that
    \begin{equation}\label{e6.27}
      R'(\varphi_\varepsilon)=R'(\varphi)+o(\varepsilon).
    \end{equation}
Finally, we have
   \begin{eqnarray}\nonumber\label{e6.28}
    K_2(\varphi_\varepsilon)&=&\frac{\partial J_1(\varphi_\varepsilon)}{\partial\eta}=\frac{n\varphi_\varepsilon^2\Big\{6[n\theta-(n-1)]\eta-2n(2\theta-1)\Big\}}{\varphi_\varepsilon^3}\\
    &&-\frac{-2n\varphi_\varepsilon\varphi_\varepsilon'\Big\{3[n\theta-(n-1)]\eta^2-2n(2\theta-1)\eta+(n\theta-1)\Big\}}{\varphi_\varepsilon^3}\\ \nonumber
    &&-\frac{n\eta(\eta-1)\Big\{[n\theta-(n-1)]\eta-(n\theta-1)\Big\}[\varphi_\varepsilon\varphi_\varepsilon''-2(\varphi_\varepsilon')^2]}{\varphi_\varepsilon^3}\\ \nonumber
    &=& K_2(\varphi)+o(\varepsilon),
   \end{eqnarray}
we \eqref{e6.26} and uniformly boundedness $\frac{\eta-1}{\varphi_\varepsilon}$ have been used again. Consequently, the mapping $T$ is continuous from $\Gamma_{\eta_0, \alpha, \beta, \gamma, \sigma}$ to $C^3([1,\eta_0])$. T show the compactness of the mapping $T$ we need only to utilize Lemma \ref{l6.8} to conclude that $T$ maps a bounded sequence $\varphi_\varepsilon$ in $\Gamma_{\eta_0, \alpha, \beta, \gamma, \sigma}$ to bounded sequences $C^4([1,\eta_0])$. So, the conclusion of the proposition follows from the Arzel\`{a}-Ascoli theorem. $\Box$\\

\noindent\textbf{Proof of Theorem \ref{t6.1}} By Proposition \ref{p6.1} and Schauder's fix point theorem, there exists a fix point $\zeta\in\Gamma_{\eta_0, \alpha, \beta, \gamma, \sigma}$ of $T$ provided $\sigma,\eta_0$ is previous chosen small, where $\alpha, \beta, \gamma$ are given in \eqref{e6.16}. By Lemma \ref{l6.2}, the parameter $\lambda(\zeta,\eta_0)=-\lambda'''$ is bounded from above and below by positive constants. Moreover, the Lipschitz continuity of third derivative of $\zeta$ follows from Lemma \ref{l6.8}. The proof was done. $\Box$\\

Finally, let's remove the dimension restriction $n\leq5$ in Theorem \ref{t6.1}. In fact, this restriction comes from proof of Lemma \ref{l6.8} about
   $$
     \frac{2}{3}+\frac{n(n-2)}{48}<1\Leftrightarrow n\leq5.
   $$
To enhance the dimension $n$, we need to compute higher derivatives of $\zeta$ as follows. At first, we rewrite \eqref{e4.7} by
  \begin{equation}\label{e6.29}
    \zeta'=L_1(\varphi)+L_2(\varphi)+\eta^2R,
  \end{equation}
where
   \begin{eqnarray*}
     L_1(\varphi)&\equiv&(\theta+1)\eta^{-1}\varphi+\Big\{[2n\theta-(2n-1)]\eta-[2n\theta-1]\Big\},\\
     L_2(\varphi)&\equiv&n\eta(\eta-1)\Big\{[n\theta-(n-1)]\eta-[n\theta-1]\Big\}\varphi^{-1},\\
     R&\equiv&\frac{\lambda_0}{\varphi}\exp\int^\eta_{\eta_0}\frac{s+1}{\varphi(s)}ds.
   \end{eqnarray*}
It' clear that $L_1(\varphi)$ is a $C^k$-function on $\eta\in[1,\eta_0]$ as long as $\varphi\in C^k([1,\eta_0])$. For any positive function $f\in C^k([1,\eta_0])$, we have $f^{-1}$ belongs also to $C^k([1,\eta_0])$. Therefore, supposing that
   $$
    f^{(l)}(1)\equiv f_l, \ \ \forall l=0,1,\cdots, k,
   $$
let's introduce a new notation for conjugate indices
   $$
    \overline{f_l}\equiv\overline{f_l}(f_0,\cdots,f_l)\equiv \frac{\partial^l f^{-1}}{\partial^l\eta}\Big|_{\eta=1}\in{\mathbb{R}}
   $$
is determined by $f_i, i=0,1,\cdots,l$ for each $l\leq k$.
Now, we can calculate the derivatives of $L_2(\varphi)$ as following.

\begin{lemm}\label{l6.9}
  If $\varphi\in C^k([1,\eta_0])$ satisfying
     \begin{equation}\label{e6.30}
      \varphi(1)=0, \ \ \varphi'(1)\not=\varphi_1\not=0, \ \ \varphi^{(l)}(1)=\varphi_l\in{\mathbb{R}}, \ \ \forall l=0,1,\cdots,k,
     \end{equation}
 the function defined by
   $$
     \phi(\eta)\equiv\begin{cases}
        \frac{\eta-1}{\varphi(\eta)}, & \eta\in(1,\eta_0]\\
         \frac{1}{\varphi_1}, & \eta=1
     \end{cases}\in C^{k-1}([1,\eta_0])
   $$
satisfies that
     \begin{equation}\label{e6.31}
        \phi_l\equiv\phi^{(l)}(1)=\overline{\Bigg(\frac{\varphi_{l+1}}{l+1}\Bigg)}, \ \ \forall l=0,1,\cdots,k-1.
     \end{equation}
\end{lemm}

\noindent\textbf{Proof.} By Taylor's expansion,
   $$
    \varphi(\eta)=\varphi_1(\eta-1)+\frac{\varphi_2}{2!}(\eta-1)^2+\cdots+\frac{\varphi_{k-1}}{(k-1)!}(\eta-1)^{k-1}+\int^\eta_1\frac{\varphi^{(k)}(t)}{(k-1)!}(\eta-t)^{k-1}dt.
   $$
Therefore,
   \begin{eqnarray*}
     f(\eta)&\equiv&\frac{\varphi(\eta)}{\eta-1}=\varphi_1+\frac{\varphi_2}{2!}(\eta-1)+\cdots+\frac{\varphi_{k-1}}{(k-1)!}(\eta-1)^{k-2}\\
       &&+\frac{\int^\eta_1\frac{\varphi^{(k)}(t)}{(k-1)!}(\eta-t)^{k-1}dt}{\eta-1}\in C^k((1,\eta_0]).
   \end{eqnarray*}
Taking derivatives $\frac{\partial^l}{\partial\eta^l}, l=1,2,\cdots,k-1$, we obtain that
   \begin{eqnarray}\nonumber\label{e6.32}
    \frac{\partial^lf}{\partial\eta^l}&=&\frac{\varphi_{l+1}}{l+1}+\frac{\varphi_{l+2}}{l+2}(\eta-1)+\cdots+\frac{\varphi_{k-1}}{(k-1)(k-2-l)!}(\eta-1)^{k-l-2}\\
     &&+\Sigma_{q=0}^l\frac{(-1)^qC_l^qq!}{(k-l+q-1)!}\frac{\int^\eta_1\varphi^{(k)}(t)(\eta-t)^{k-l+q-1}dt}{(\eta-1)^{q+1}}, \ \ \forall l=1,2,\cdots,k-2
   \end{eqnarray}
and
    \begin{equation}\label{e6.33}
      \frac{\partial^{k-1}f}{\partial\eta^{k-1}}=\Sigma_{q=0}^{k-1}(-1)^qC_{k-1}^q\frac{\int^\eta_1\varphi^{(k)}(t)(\eta-t)^qdt}{(\eta-1)^{q+1}}.
    \end{equation}
Noting that for $l=1,2,\cdots,k-2$,
    $$
     \lim_{\eta\to1^+}\frac{\partial^lf}{\partial\eta^l}(\eta)=\frac{\varphi_{l+1}}{l+1},
    $$
we have $f\in C^{k-2}([1,\eta_0])$ and
   \begin{equation}\label{e6.34}
     f_l\equiv\frac{\partial^lf}{\partial\eta^l}(1)=\frac{\varphi_{l+1}}{l+1}, \ \ \forall l=0,1,\cdots,k-2
   \end{equation}
by Cauchy's theorem. To show $f\in C^{k-1}([1,\eta_0])$, one need only use \eqref{e6.33} to deduce that
   $$
    \lim_{\eta\to1^+}\frac{\partial^{k-1}f}{\partial\eta^{k-1}}(\eta)=\lim_{\xi\to1^+}\varphi^{(k)}(\xi)\Sigma_{q=0}^{k-1}(-1)^q\frac{C_{k-1}^q}{q+1}=\frac{\varphi_k}{k}
   $$
by Cauchy's theorem, and thus conclude that $f\in C^{k-1}([1,\eta_0])$ and
   \begin{equation}\label{e6.35}
     f_{k-1}\equiv\frac{\partial^{k-1}f}{\partial\eta^{k-1}}(1)=\frac{\varphi_k}{k}.
   \end{equation}
 So, \eqref{e6.31} follows from conjugatings of \eqref{e6.34} and \eqref{e6.35}. $\Box$\\

Now, supposing that
   \begin{equation}\label{e6.36}
    \zeta^{(k+1)}=\frac{\partial^kL_1(\varphi)}{\partial\eta^k}+\frac{\partial^kL_2(\varphi)}{\partial\eta^k}+R\cdot R_{k+1},
   \end{equation}
one has the induction formula
   \begin{equation}\label{e6.37}
     \begin{cases}
      \displaystyle R_{k+1}=\Bigg(-\frac{\varphi'}{\varphi}+\frac{\eta+1}{\varphi}\Bigg)+\frac{\partial R_k}{\partial\eta}, & \forall k\geq1,\\
       R_1=\eta^2.
     \end{cases}
   \end{equation}
Setting
   $$
    F(\eta)\equiv-\frac{\varphi'}{\varphi}+\frac{\eta+1}{\varphi}
   $$
for simplicity, we have the following result by Lemma \ref{l6.9}.

\begin{lemm}\label{l6.10}
  Letting $\varphi^{k+1}([1,\eta_0])$ satisfying \eqref{e6.30}, one has
    \begin{equation}\label{e6.38}
      F_l\equiv\frac{\partial^lF}{\partial\eta^l}(1)=\Sigma_{q=0}^lC_l^q\overline{\Bigg(\frac{\varphi_{q+1}}{q+1}\Bigg)}\Bigg(1-\frac{\varphi_{l-q+2}}{l-q+1}\Bigg), \ \ \forall l=0,1,\cdots,k-1.
    \end{equation}
\end{lemm}

\noindent\textbf{Proof.} Rewritten $F$ by
   $$
    F(\eta)=\frac{\eta-1}{\varphi}\Bigg(1-\frac{\varphi'-2}{\eta-1}\Bigg),
   $$
it yields from \eqref{e6.31}, \eqref{e6.34} and \eqref{e6.35} in Lemma \ref{l6.9} that \eqref{e6.38} holds. $\Box$\\

As a corollary of Lemma \ref{l6.9} and \ref{l6.10}, one has
   \begin{eqnarray}\nonumber\label{e6.39}
     R_{k+1}\Big|_{\eta=1}&=&F_0+F_1+\cdots+F_{k-1}\\
      &=&\psi_1(\varphi_0,\varphi_1,\cdots,\varphi_k)-\frac{\varphi_{k+1}}{2k}
   \end{eqnarray}
when $k\geq3$, where $\psi_1(\varphi_0,\varphi_1,\cdots,\varphi_k)$ is a bounded function depending only on $\varphi_0,\varphi_1,\cdots,\varphi_k$. Furthermore,
  \begin{eqnarray}\label{e6.40}
     \frac{\partial^kL_1(\varphi)}{\partial\eta^k}\Big|_{\eta=1}&=&\psi_2(\varphi_0,\varphi_1,\cdots,\varphi_k)\\ \nonumber
     \frac{\partial^kL_2(\varphi)}{\partial\eta^k}\Big|_{\eta=1}&=&\psi_3(\varphi_0,\varphi_1,\cdots,\varphi_k)-n(n-2)\overline{\Bigg(\frac{\varphi_{k+1}}{k+1}\Bigg)}
  \end{eqnarray}

\begin{lemm}\label{l6.11}
  For any positive function $f\in C^k([1,\eta_0])$, one has
    \begin{equation}\label{e6.41}
      \overline{f_k}\equiv\frac{\partial^kf^{-1}}{\partial\eta^k}\Big|_{\eta=1}=-\frac{f_k}{f_0^2}.
    \end{equation}
\end{lemm}

Substituting \eqref{e6.39}-\eqref{e6.41} into \eqref{e6.36}, one concludes that
   \begin{equation}\label{e6.42}
     \lim_{\eta\to1^+}\zeta^{(k+1)}(\eta)=\psi(\varphi_0,\varphi_1,\cdots,\varphi_k)+\Bigg[-\frac{1}{2k}\lim_{\eta\to1^+}R+\frac{n(n-2)}{\varphi_1^2(k+1)}\Bigg]\varphi^{(k+1)}(1).
   \end{equation}
Now, for any given $n\geq10$, let's take a positive integer $k$ so large that
   \begin{equation}\label{e6.43}
     \varsigma\equiv\frac{1}{2k}\Bigg[4+\frac{n(n-2)}{2}\Bigg]-\frac{n(n-2)}{4(k+1)}\in(0,1),
   \end{equation}
and introduce a family of functions
    \begin{eqnarray*}
     &\Gamma_{\eta_0, \alpha, \gamma, \sigma}\equiv\Big\{\varphi\in C^{k+1}([1,\eta_0])\Big|\ \varphi^{(l)}(1)=\alpha_l, \ \ \forall l=0,1,\cdots,k+1&\\
     &\varphi^{(l)}(\eta)\in[\alpha_l-\sigma,\alpha_l+\sigma], \ \ \forall \eta\in(1,\eta_0], \ \ l=0,1,\cdots,k&\\
     &\varphi^{(k+1)}(\eta)\in[\gamma-1,\gamma+1]\Big\}&
    \end{eqnarray*}
for any $\alpha=(\alpha_0,\alpha_1,\cdots,\alpha_{k+1})\in{\mathbb{R}}^{k+1}, \eta_0>1, \gamma,\sigma>0$, it's clear that closed convex subset of $C^{k+1}([1,\eta_0])$ endowed with norm $||\cdot||_{C^{k+1}([1,\eta_0])}$. A similar argument as in proof of Theorem \ref{t6.1} and a application of \eqref{e6.42} show that for
  \begin{equation}\label{e6.44}
   \begin{cases}
    \alpha_0=0, \ \ \alpha_1=2, \\[5pt] \displaystyle\alpha_l=\frac{4l(l-1)\psi(\alpha_0,\alpha_1,\cdots,\alpha_{l-1})}{4l(l+1)+n(n-2)}, & \forall l=3,\cdots,k+1
   \end{cases}
  \end{equation}
and
   \begin{equation}\label{e6.45}
     \gamma=\frac{4k(k+1)\psi(\alpha_0,\cdots,\alpha_k)}{4(k+1)(k+2)+n(n-2)},
   \end{equation}
the mapping $T$ defined by \eqref{e6.2} and \eqref{e6.3} is a continuous and compact mapping from $\Gamma_{\eta_0, \alpha, \gamma, \sigma}$ into $\Gamma_{\eta_0, \alpha, \gamma, \sigma}$.  After applying Schauder's fix point theorem, we arrive at the following theorem.

\begin{theo}\label{t6.2}
  For any $n\geq2$ and $k\geq3$ so large that \eqref{e6.43} holds, we assume the parameters $\alpha, \gamma$ are given by \eqref{e6.44} and \eqref{e6.45}. If $\eta_0-1$ is small, there exists solution $\zeta\in C^{k+1}([1,\eta_0])$ of \eqref{e4.7} for some negative number $\lambda_0$.
\end{theo}

  \vspace{40pt}

\section{Long time existence and quadratic growth}

 In this section, let's extend $\zeta$ to be a global solution as follows.

 \begin{theo}\label{t7.1}
   The local positive solution in Theorem \ref{t6.2} for $n\geq2, \theta>(n-7)/n^2, k\geq3$ can be extended to be a global positive solution $\zeta\in C^{k+1}([1,+\infty))\cap C^\infty((1,+\infty))$. Moreover, there exist a large constant $\eta_1(\theta)$ and a small constant $\varepsilon_0(\eta_0,\eta_1)>0$ such that
     \begin{equation}\label{e7.1}
       \zeta(\eta)>\varepsilon_0\eta^2, \ \ \forall \eta\geq\eta_1.
     \end{equation}
 \end{theo}

 \noindent\textbf{Proof.} At first, it's clear that the solution $\zeta$ can be extended beyond fixed point $\eta=\eta_1>\eta_0$ supposing that it is bounded from above and below by positive constants. Secondly, we claim that if one choose $\varrho=\varrho(n)>0$ small, there holds
    \begin{equation}\label{e7.2}
      \zeta(\eta)\geq\varrho(\eta-1), \ \ \forall \eta\geq1
    \end{equation}
  for fixed $\eta_0-1$ small. To show \eqref{e7.2}, one needs only to prove that $\xi(\eta)=\varrho(\eta-1)$ is a sub-solution to \eqref{e4.7} for $\lambda'''=-\lambda(\zeta,\eta_0)\equiv\lambda_0$. In fact,
   \begin{eqnarray*}
     && L\xi\equiv(\theta+1)\eta^{-1}\xi+n\eta(\eta-1)\Big\{[n\theta-(n-1)]\eta-[n\theta-1]\Big\}\xi^{-1}\\
     &&+\Big\{[2n\theta-(2n-1)]\eta-[2n\theta-1]\Big\}+\frac{\lambda_0\eta^2}{\xi}\exp\int^\eta_{\eta_0}\frac{s+1}{\xi(s)}ds\\
     &>&\frac{n}{\varrho}\eta\Big\{[n\theta-(n-1)]\eta-[n\theta-1]\Big\}+\Big\{[2n\theta-(2n-1)]\eta-[2n\theta-1]\Big\}\\
     &&+\frac{[32+4n(n-2)]\eta^2}{\varrho[(\eta_0-1)(\alpha+\sigma)+4]}e^{\frac{\eta-\eta_0}{b}}, \forall \eta>1
   \end{eqnarray*}
 by Lemma \ref{l6.2}. Since for small $\alpha,\sigma$ and then small $\eta_0-1$,
    $$
     \frac{32+4n(n-2)}{(\eta_0-1)(\alpha+\sigma)+4}\geq7+n(n-2),
    $$
 we have
    \begin{eqnarray}\nonumber\label{e7.3}
      L\xi&>&\frac{n}{\varrho}\eta\Big\{[n\theta-(n-1)]\eta-[n\theta-1]\Big\}+\Big\{[2n\theta-(2n-1)]\eta-[2n\theta-1]\Big\}\\
       &&+\frac{7+n(n-2)}{\varrho}\eta^2\equiv f(\eta), \ \ \forall \eta>1.
    \end{eqnarray}
 Thus, using the fact of positivity of quadratic function $f(\eta)-\varrho$ for all $\eta\geq1$ when $\theta>(n-7)/n^2$ and $\varrho$ small, we get
    \begin{equation}\label{e7.4}
     L\xi\geq \varrho=\xi', \ \ \forall\eta>1.
    \end{equation}
 Noting that
    \begin{equation}\label{e7.5}
       \zeta(\eta_0)>\xi(\eta_0)
    \end{equation}
holds for any
   $$
     \zeta\in\Gamma_{\eta_0, \alpha, \beta, \gamma, \sigma},
   $$
one obtain a a-priori positive lower $\xi$ for $\zeta$ by comparison of $\zeta$ and $\xi$ for $\eta>\eta_0$. Thirdly, for any $\eta<\eta_2, \eta_2>\eta_0$, it's inferred from the boundedness of $\zeta$ from below that
   $$
    \zeta'\leq(\theta+1)\zeta+C_{\eta_1}.
   $$
 Thus, $\zeta$ has also a-priori bounded from above by positive functions. Namely, any local positive solution in Theorem \ref{t6.2} can be extended to be a global positive solution on $[1,+\infty)$ and thus belongs to $C^\infty((1,+\infty))$, as long as $\theta>(n-7)/n^2$. To show \eqref{e7.1} holds, let's prove a lemmas.

 \begin{lemm}\label{l7.1}
   Under the assumption of Theorem \ref{t7.1}, there exist a large constant $\eta_1(\theta)$ and a small constant $\varepsilon_0(\eta_0,\eta_1)>0$ such that \eqref{e7.1} holds.
 \end{lemm}

 \noindent\textbf{Proof.} At first, one choose $\eta_1$ so large that
   $$
   \begin{cases}
     n(\eta-1)\Big\{[n\theta-(n-1)]\eta-[n\theta-1]\Big\}\geq n[n\theta-n]\eta^2, & \forall\eta\geq\eta_1\\
     [2n\theta-(2n-1)]\eta-[2n\theta-1]\geq[2n\theta-2n]\eta, & \forall\eta\geq\eta_1.
   \end{cases}
   $$
 Thus, the function $\xi(\eta)\equiv\varepsilon_0\eta^2$ satisfies that
   \begin{eqnarray*}
     L\xi&=&(\theta+1)\varepsilon_0\eta+n[n\theta-n]\eta+[2n\theta-2n]\eta\\
      &&+\frac{\lambda_0}{\varepsilon_0}\exp\Bigg\{\frac{1}{\varepsilon_0}\int^\eta_{\eta_0}\frac{1}{s}+\frac{1}{s^2}ds\Bigg\}\\
      &\geq&\Bigg\{(\theta+1)\varepsilon_0+(n^2+2n)(\theta-1)+\frac{7(\eta_0-1)}{\varepsilon_0\eta_0}\Bigg\}\eta\\
      &\geq&2\varepsilon_0\eta=\xi'
   \end{eqnarray*}
 once $\varepsilon_0$ is chosen small, where Lemma \ref{l6.2} has been used again. Shrinking the constant $\varepsilon_0$ again with respect to $\eta_1$, one may assume that
    $$
     \xi(\eta_1)<\zeta(\eta_1).
    $$
 As a result, we obtain \eqref{e7.1} by comparing $\zeta$ with $\xi$. The proof was done. $\Box$\\

  As an application of Theorem \ref{t7.1}, we have the following result of finite time blowup.

\begin{coro}\label{c7.1}
   Under the assumptions of Theorem \ref{t7.1}, the solution $\zeta$ given by Theorem \ref{t6.2} and \ref{t7.1} yields a positive solution $\eta(t)$ of \eqref{a4.6} which blows up at finite time $t=T_0\in{\mathbb{R}}$. More precisely,
     \begin{equation}\label{e7.6}
       T_\infty\equiv\int^{+\infty}_{\eta_0}\frac{ds}{\zeta(s)}<+\infty.
     \end{equation}
\end{coro}

Letting $\zeta_{\eta_0}$ be the solution obtained in Theorem \ref{t6.2} and \ref{t7.1} with respect to $\eta_0$, we have $\lambda(\zeta_{\eta_0},\eta_0)$ is small as long as $\eta_0-1$ is small. So, we have the following upper bound of $\zeta$.

\begin{theo}\label{t7.2}
  Under the assumptions of Theorem \ref{t7.1}, for any $\theta\in[1/n,n/(n+1))$, if one choose $\eta_0-1$ small enough, there holds
     \begin{equation}\label{e7.7}
        \zeta(\eta)\leq\eta^2, \ \ \forall \eta\geq\eta_2
     \end{equation}
  for some large number $\eta_2(n,\theta)$.
\end{theo}

Before proving the theorem, let's recall from Lemma \ref{l6.2} that there exists a monotone decreasing sequence of $\eta_0(l)>1, l\in{\mathbb{N}}$ satisfying
    \begin{equation}\label{e7.8}
     \eta_0(l)\downarrow1, \ \ \mbox{ as } l\uparrow+\infty
    \end{equation}
such that
    \begin{equation}\label{e7.9}
      \lambda(l)\equiv\lambda(\zeta_{\eta_0(l)},\eta_0(l))\downarrow 0, \ \ \mbox{ as } l\uparrow+\infty.
    \end{equation}
Utilizing the comparison principle for \eqref{e4.7}, one has also
    \begin{equation}\label{e7.10}
       \zeta_l\equiv\zeta_{\eta_0(l)}\downarrow\ \ \mbox{as } l\uparrow+\infty.
    \end{equation}
Now, we impose first the following lemma.

\begin{lemm}\label{l7.2}
  Under the assumptions of Theorem \ref{t7.2}, for any $k\in{\mathbb{N}}$, there exist $\eta_k>k$ and $l=l_k$, such that
     \begin{equation}\label{e7.11}
       \zeta(\eta)<\eta^2, \ \ \eta=\eta_k, \ \ l=l_k.
     \end{equation}
\end{lemm}

\noindent\textbf{Proof.} Suppose on the contrary, there exists $k_*\in{\mathbb{N}}$ such that
   \begin{equation}\label{e7.12}
     \zeta_l(\eta)\geq\eta^2, \ \ \forall \eta\geq k_*
   \end{equation}
holds for all $l\in{\mathbb{N}}$. By \eqref{e4.7} for $\lambda'''=-\lambda_l$ and \eqref{e7.2} for $\varrho(n)$, we have
   \begin{eqnarray*}
     \zeta_l'&\leq&(\theta+1)\eta^{-1}\zeta_l+n[n\theta-(n-1)]_+(\eta-1)+[2n\theta-(2n-1)]_+\eta\\
      &&+\lambda_l\exp\Bigg\{\int^{k_*}_{\eta_0}\frac{s+1}{\varrho(n)(s-1)}ds+\int^\eta_{k_*}\frac{s+1}{s^2}ds\Bigg\}\\
      &\leq&(\theta+1)\eta^{-1}\zeta_l+\Bigg\{n[n\theta-(n-1)]_++[2n\theta-(2n-1)]_+\Bigg\}\eta\\
      &&+\frac{\lambda_l\eta}{\eta_*}\Bigg(\frac{k_*}{\eta_0}\Bigg)^2\exp\Bigg\{\frac{1}{\varrho(n)}(k_*-\eta_0)+\frac{1}{k_*}\Bigg\}, \ \ \forall\eta\geq k_*
   \end{eqnarray*}
holds for all $l\in{\mathbb{N}}$. Using the fact of \eqref{e7.9} and $\theta<n/(n+1)$, one has
   $$
     \zeta_l'(\eta)<(\theta+1)\eta^{-1}\zeta+(1-\theta-\varrho_*)\eta, \ \ \forall\eta>k_*
   $$
for some positive constant $\varrho_*$. As a result, one concludes that
   $$
     \zeta(\eta)\leq\frac{1-\theta-\varrho_*}{1-\theta}\eta^2+C_*\eta^{1+\theta}, \ \ \forall\eta\geq k_*
   $$
for some positive constant $C_*$. This contradicts with \eqref{e7.12} for large $\eta$. The conclusion of the lemma was drawn. $\Box$\\

Now, one can complete the proof of Theorem \ref{t7.2}.\\

\noindent\textbf{Continue the proof of Theorem \ref{t7.2}.} Noting first the integration O.D.E. has no comparison lemma to be used, thus we first change the equation into second order form without integration term as following
   \begin{eqnarray}\nonumber\label{a.1}
     &\log L_1(\zeta)=\log\lambda_0+2\log\eta+\int^\eta_{\eta_0}\frac{s+1}{\zeta(s)}ds&\\
     &\Leftrightarrow\zeta\zeta''+(\zeta')^2-\Big\{\Big[2n\theta-(2n-1)\Big]\eta-[2n\theta-1]\Big\}\zeta'&\\ \nonumber
     &+L_0(\zeta)-g'(\eta)=\Bigg(\frac{2}{\eta}+\frac{\eta+1}{\zeta(\eta)}\Bigg)L_1(\zeta),&
   \end{eqnarray}
where $\zeta\equiv\zeta_l$,
    $$
     g(\eta)\equiv\Bigg\{n\eta(\eta-1)\Big\{\Big[n\theta-(n-1)\Big]\eta-[n\theta-1]\Big\}\Bigg\}
    $$
and
    \begin{eqnarray*}
      L_0(\zeta)&\equiv&(\theta+1)\eta^{-2}\zeta^2-2(\theta+1)\eta^{-1}\zeta\zeta'-\Big[2n\theta-(2n-1)\Big]\zeta\\
      L_1(\zeta)&\equiv&\zeta\zeta'-(\theta+1)\eta^{-1}\zeta^2-\zeta\Big\{\Big[2n\theta-(2n-1)\Big]\eta-[2n\theta-1]\Big\}-g(\eta).
    \end{eqnarray*}
Now, we use maximum type principle to deduce that
   \begin{equation}\label{a.2}
     \zeta(\eta)\leq\eta^2, \ \ \forall\eta\in[\eta_k,\eta_{k+1}]
   \end{equation}
for $k$ large. In fact, assume on the contrary, there exists $\eta_*\in(\eta_k,\eta_{k+1})$ such that
   $$
    (\zeta-\eta^2)\Big|_{\eta=\eta_*}=\max_{\eta\in[\eta_k,\eta_{k+1}]}(\zeta-\eta^2)\geq0.
   $$
Therefore,
   \begin{equation}\label{a.3}
      \zeta(\eta_*)\geq\eta_*^2, \ \ \zeta'(\eta_*)=2\eta_*.
   \end{equation}
To conclude the possibility
   \begin{equation}\label{a.4}
     \zeta(\eta_*)>\eta_*^2,
   \end{equation}
we need the following two lemmas. The first told us that $\eta^2$ is a super-solution of \eqref{a.1}.

\begin{lemm}\label{l7.3}
  Under the assumption of Theorem \ref{t7.2}, $\eta^2$ satisfies that
    \begin{equation}\label{a.5}
     (L.H.S.-R.H.S)(\eta^2)=-(n-1)\big[(n+1)\theta+1\big]\eta+o(\eta)<0
    \end{equation}
  if $\theta>-\frac{1}{n+1}$ and $\eta$ large.
\end{lemm}

\noindent\textbf{Proof.} Substituting $\eta^2$ into both sides of \eqref{a.1}, one gets that
  \begin{eqnarray*}
    L.H.S.&\sim& 2\eta^2+4\eta^2+(\theta+1)\eta^2-4(\theta+1)\eta^2-\Big[2n\theta-(2n-1)\Big]\eta^2\\
    &&-2\Big[2n\theta-(2n-1)\Big]\eta^2-3n\Big[n\theta-(n-1)\Big]\eta^2\\
    &&+\Big[2(2n\theta-1)+2n(2n\theta-n)\Big]\eta\\
    R.H.S.&\sim&\frac{3}{\eta}\Bigg\{2\eta^3-(\theta+1)\eta^3-\Big[2n\theta-(2n-1)\Big]\eta^3-n\Big[n\theta-(n-1)\Big]\eta^3\\
      &&-3\big[(2n\theta-1)+n(2n\theta-n)\big]\eta-\big[-(n+1)^2\theta+n(n+1)\big]\eta.
  \end{eqnarray*}
So, \eqref{a.5} holds. $\Box$\\

\noindent\textbf{Continue the proof of Theorem \ref{t7.2}} Setting $\xi\equiv\eta^4-\zeta^2$, we will prove that $\xi$ is non-positive on $[\eta_k,\eta_{k+1}]$ for $k$ large. If not, by shrinking the interval, one may assume $\xi$ is positive inside $(\eta_k,\eta_{k+1})$ and attains a positive maximum at $\eta_*$. By \eqref{a.1}, there holds
   \begin{equation}\label{a.6}
     \frac{1}{2}\xi''+a_1(\eta)\xi'+a_2(\eta)\xi+b(\eta)-\Bigg(-\frac{\eta+1}{\zeta}+\frac{\eta+1}{\eta^2}\Bigg)L_1(\eta^2)=0,
   \end{equation}
where
  \begin{eqnarray*}
    a_1(\eta)&\equiv&-\frac{\big[2n\theta-(2n-1)\big]\eta-\big[2n\theta-1\big]}{\zeta}-\Bigg(\frac{1}{\eta}+\frac{\eta+1}{2\zeta}\Bigg)\\
    a_2(\eta)&\equiv&(\theta+1)\eta^{-2}-\frac{2n\theta-(2n-1)}{\zeta+\eta^2}+\Bigg(\frac{2(\theta+1)}{\eta^2}+\frac{(\theta+1)(\eta+1)}{\zeta\eta}\Bigg)\\
    &&+\Bigg(\frac{2}{\eta}+\frac{\eta+1}{\zeta}\Bigg)\frac{\big[2n\theta-(2n-1)\big]\eta-(2n\theta-1)}{\zeta+\eta^2}
  \end{eqnarray*}
are both small function as long as $\eta_k$ large, and
  $$
   b(\eta)\equiv-4\eta^3\Bigg(\frac{1}{\xi}-\frac{1}{\eta^2}\Bigg)\Bigg\{\Big[2n\theta-(2n-1)\Big]\eta-\Big[2n\theta-1\Big]\Bigg\}
  $$
is negative as long as $\theta<\frac{2n-1}{2n}$ and $\eta$ large. Noting that for $\theta<\frac{n}{n+1}$, there holds
   \begin{eqnarray*}
     L_1(\eta^2)=(n+1)\Big[-(n+1)\theta+n\Big]>0,
   \end{eqnarray*}
one concludes from \eqref{a.6} that
   \begin{equation}\label{a.7}
     \frac{1}{2}\xi''+a_1(\eta)\xi'+a_2(\eta)\xi\geq0, \ \ \forall \eta\in(\eta_k,\eta_{k+1})
   \end{equation}
and
   \begin{equation}\label{a.8}
     \xi(\eta_k)=\xi(\eta_{k+1})=0.
   \end{equation}
By smallness of $a_1$ and $a_2$ for large $k$, the largest eigenvalue of $\frac{1}{2}\frac{d^2}{d\eta^2}+a_1(\eta)\frac{d}{d\eta}+a_2(\eta)$ is negative for Dirichlet problem, it yields from \eqref{a.7} and \eqref{a.8} that $\xi$ is non-positive everywhere inside $[\eta_k,\eta_{k+1}]$, contradiction holds. Theorem \ref{t7.2} holds true. $\Box$\\

\begin{coro}\label{c7.2}
  Under the assumptions of Theorem \ref{t7.2}, the restored solution $u$ of \eqref{e4.4} satisfies a large condition
    \begin{equation}\label{e7.10}
      \lim_{r\to R_\infty^-}u(r)=+\infty
    \end{equation}
  on boundary, where $R_\infty\equiv e^{T_\infty}<+\infty$.
\end{coro}

\noindent\textbf{Proof.} By \eqref{e7.7} and \eqref{e5.3}, one has
  $$
   t= T_\infty-\int^{+\infty}_\eta\frac{ds}{\zeta(s)}\leq T_\infty-\frac{1}{\eta}
  $$
for $\eta$ large. Therefore, there holds
   $$
    \eta(t)\geq\frac{1}{T_\infty-t}, \ \ \forall t\in(T_*,T_\infty)
   $$
or equivalent
    $$
     (\log v)_r=w\geq\frac{1}{r(T_\infty-\log r)}, \ \ \forall r\in(e^{T_*},R_\infty).
    $$
Thus,
   \begin{eqnarray*}
    \frac{v}{v_0}&\geq&\exp\int^r_{r_0}\frac{1}{r(T_\infty-\log r)}ds\\
     &=&\frac{T_\infty-\log r_0}{T_\infty-\log r}, \ \ \forall r\in(e^{T_*},R_\infty).
   \end{eqnarray*}
Integrating over $r$, one gets \eqref{e7.10} thanks to
   $$
    \int^{R_\infty}_{R_\infty/2}\frac{T_\infty-\log r_0}{T_\infty-\log r}dr=+\infty.
   $$
$\Box$\\

\vspace{40pt}

\section{False of Bernstein property under Euclidean completeness}

Letting $v$ be the function given by \eqref{e4.10}, we have
   \begin{eqnarray*}
     \int^{v_0}_{v(r)}\frac{dv}{\sqrt{v^3}}<\sqrt{a}r\Rightarrow v(r)>\frac{4}{\Big(\frac{2}{\sqrt{v_0}}+\sqrt{a}r\Big)^2}, \ \ \forall r\geq0.
   \end{eqnarray*}
Twice integrating over $r$, the solution
  $$
    u(r)=\int^r_0\int^t_0v(s)dsdt=\int^r_0(r-s)v(s)ds\geq4\int^r_0\frac{r-s}{\Big(\frac{2}{\sqrt{v_0}}+\sqrt{a}s\Big)^2}ds
  $$
satisfies the large condition
  $$
    \lim_{r\to+\infty}u(r)=+\infty
  $$
when $\theta>1/2$. Equivalently, the solution $\varphi(x)\equiv u(|x|)$ of \eqref{e4.3} for positive pair $(\lambda,1)$ satisfies a large condition
  \begin{equation}\label{e8.1}
    \lim_{x\to\infty}\varphi(x)=+\infty.
  \end{equation}

By Proposition \ref{p5.2}, the solution $\psi(y)\equiv u(|x|)$ obtained in Theorem \ref{t6.2} and \ref{t7.1} belongs to $C^5(B_{R_\infty})$, where $R_\infty$ is given in \eqref{e7.9}. Furthermore, Corollary \ref{c7.2} tells us that $\psi(y)$ satisfies a large condition
  \begin{equation}\label{e8.2}
    \lim_{y\to\partial B_{R_\infty}}\psi(y)=+\infty
  \end{equation}
provided $\theta\in[1/n,n/(n+1))$. Consequently, as illustrated in Section 4, after scaling $\varphi, \psi$ such that they are solutions to opposite pairs, the function
  $$
   u(x,y)\equiv\varphi(x)+\psi(y), x\in{\mathbb{R}}, y\in B_{R_\infty}\subset{\mathbb{R}}^2
  $$
is a $C^5(\Omega)$ solution to \eqref{e1.1} on convex domain
   $$
    \Omega\equiv{\mathbb{R}}\times B_{R_\infty}.
   $$
It follows from the large condition
   $$
    \lim_{(x,y)\to\partial\Omega}u(x,y)=+\infty
   $$
that the graph of $u$ is a Euclidean complete hypersurface to \eqref{e1.1} for $\theta\in(1/2,n/(n+1))$. Therefore, the following counter example to Bernstein problem holds.

\begin{theo}\label{t8.1}
  For $N\geq3$ and
   \begin{equation}\label{e8.3}
       \theta\in(1/2,(N-1)/N),
    \end{equation}
   there exists a $C^5$ Euclidean complete hypersurface on ${\mathbb{R}}^{N+1}$ which is not a elliptic paraboloid.
\end{theo}

\noindent\textbf{Proof.} Above discussions actually give the counter example for $N=n+1\geq3$ and $\theta\in(1/2,(N-1)/N)$. $\Box$\\

\noindent\textbf{Remark 8.1} As shown by Trudinger-Wang in \cite{TW2}, affine complete locally convex solution to \eqref{e1.1} must be Euclidean complete. Therefore, a corollary of Theorem \ref{t8.1} is that for $N\geq3$ and $\theta$ satisfying \eqref{e8.3}, Bernstein theorem under affine complete fails to hold.

\vspace{40pt}

\section*{Acknowledgments}

The author would like to express his deepest gratitude to Professors Xi-Ping Zhu, Kai-Seng Chou, Xu-Jia Wang and Neil Trudinger for their constant encouragements and warm-hearted helps. This paper is also dedicated to the memory of Professor Dong-Gao Deng.\\

\end{document}